\documentclass{amsart}

\usepackage{epsfig,amsmath}
\usepackage{xspace}
\usepackage[psamsfonts]{amssymb}
\usepackage[latin1]{inputenc}
\usepackage{color}

\usepackage[normalem]{ulem} % pour la commande \sout

\usepackage{hyperref}

\newtheorem{proposition}{\bf Proposition}
\newtheorem{definition}{\bf Definition}
\newtheorem{theorem}{\bf Theorem}
\newtheorem{lemma}{\bf Lemma}
\newtheorem{corollary}{\bf Corollary}
\newtheorem{conjecture}{\bf Conjecture}

\definecolor{dkGreen}{rgb}{0,0.5,0}

\renewcommand{\epsilon}{\varepsilon} % à éviter : préférer \e
\newcommand{\egaldef}{\overset{\text{\tiny def}}{=}}
\def\tend{\longrightarrow}

\def\ds{\displaystyle}
\def\wt{\widetilde}
\def\on{\operatorname}

\def\C{{\mathbb C}}
\def\D{{\mathbb D}}
\def\H{{\mathbb H}}
\def\Q{{\mathbb Q}}
\def\Z{{\mathbb Z}}
\def\R{{\mathbb R}}
\def\P{{\mathbb P}^1}
\def\N{{\mathbb N}}
\def\U{{\mathbb U}}

\def\d{\delta}
\def\a{\alpha}
\def\e{\varepsilon}
\def\cal{\mathcal}
\def\re{{\on{Re}}}
\def\im{{\on{Im}}}
\def\rad{{\on{rad}}}
\def\Re{{\on{Re}}}
\def\Im{{\on{Im}}}
\def\trunc{{\on{trunc}}}
\def\Cremer{Cremer\xspace}

\def\proof{\noindent{\bf Proof. }}
\def\remark{\vskip.2cm \noindent{\bf Remark. }}
\def\qedprop{\hfill{\hbox{%
  \hskip 1pt%
  \vrule width 7pt height 6pt depth 1.5pt%
  \hskip 1pt}} \vskip.3cm}
\def\qedlem{\hfill{$\square$} \vskip.3cm}

\def\ds{\displaystyle}

\begin{document}

\title[The Brjuno Function and
the Size of Siegel Disks.]{The Brjuno Function Continuously Estimates
the Size of Quadratic Siegel Disks.}
\subjclass{}
\begin{author}[X.~Buff]{Xavier Buff}
\email{buff$@$picard.ups-tlse.fr}
\address{ %
  Universit\'e Paul Sabatier\\
  Laboratoire Emile Picard \\
  118, route de Narbonne \\
  31062 Toulouse Cedex \\
  France }
\end{author}
\begin{author}[A.~Chéritat]{Arnaud Chéritat}
\email{cheritat$@$picard.ups-tlse.fr}
\address{ %
  Universit\'e Paul Sabatier\\
  Laboratoire Emile Picard \\
  118, route de Narbonne \\
  31062 Toulouse Cedex \\
  France }
\end{author}

\begin{abstract}
If $\alpha$ is an irrational number, we define Yoccoz's Brjuno function $\Phi$
by
$$\Phi(\a)=\sum_{n\geq 0} \a_0\a_1\cdots\a_{n-1}\log\frac{1}{\a_n},$$
where $\a_0$ is the fractional part of $\a$ and $\a_{n+1}$ is the fractional
part of ${1/\a_n}$. The numbers $\a$ such that $\Phi(\a)<\infty$ are called
the Brjuno numbers.

The quadratic polynomial $P_\alpha:z\mapsto e^{2i\pi \alpha}z+z^2$ has an
indifferent fixed
point at the origin. If $P_\alpha$ is linearizable, we let $r(\alpha)$ be the
conformal radius of the Siegel disk and we set $r(\alpha)=0$ otherwise.

Yoccoz \cite{y} proved that $\Phi(\a)=\infty$ if and only if
$r(\a)=0$ and that the restriction of $\a\mapsto \Phi(\a)+\log
r(\a)$ to the set of Brjuno numbers is bounded from below by a
universal constant. In \cite{bc2}, we proved that it is also
bounded from above by a universal constant. In fact, Marmi, Moussa
and Yoccoz \cite{mmy} conjecture that this function extends to
$\R$ as a H\"older function of exponent $1/2$. In this article, we
prove that there is a continuous extension to $\R$.
\end{abstract}

\maketitle

\section{Introduction.}

For any irrational number $\alpha\in \R\setminus \Q$, we denote by
$(p_n/q_n)_{n\geq 0}$ the approximants to $\alpha$ given by its
continued fraction expansion (by convention, $p_0=\lfloor
\a\rfloor$ is the integer part of $\a$ and $q_0=1$).

\remark Every time we use the notation $p/q$ for a rational
number, we mean that $q>0$ and $p$ and $q$ are coprime. \vskip.2cm

We denote by $\lfloor\alpha\rfloor\in \Z$ the integer part of
$\alpha$, i.e., the largest integer $n \leq \alpha$, by
$\{\alpha\} = \alpha - \lfloor\alpha\rfloor$ the fractional part
of $\alpha$, and we define $( \alpha_n)_{n\geq 0}$ recursively by
setting $\alpha_{0} = \{\alpha\}$ and $\alpha_{n+1}=\{1/\alpha_n\}.$
We then define $
\beta_{-1}=1$ and $ \beta_{n}=
 \alpha_0 \alpha_{1}\cdots  \alpha_{n}$.

\begin{definition}{\sc (The Yoccoz function).}
If $\alpha$ is an irrational number, we define
$$\Phi(\alpha) = \sum_{n=0}^{+\infty} \beta_{n-1}
\log \frac{1}{\alpha_{n}}.$$ If $\alpha$ is a rational number we
define $\Phi(\alpha)=+\infty$. Irrational numbers for which
$\Phi(\alpha)<\infty$ are called Brjuno numbers. Other irrational
numbers are called \Cremer numbers.
\end{definition}

\remark The set ${\cal B}$ of Brjuno numbers has full measure in
$\R$. It contains the set of all Diophantine numbers, i.e.,
numbers for which $\log q_{n+1} = {\cal O}(\log q_n).$ \vskip.2cm

%The reason why the Yoccoz function is of interest is that it gives
%an estimate on the size of the Siegel disks of
We study the quadratic polynomials
\[P_\a:z\mapsto e^{2i\pi \a}z+z^2\]
for $\a\in \R$.
It is known that such $P_\alpha$ is linearizable -- and so, has a
Siegel disk -- if and only if $\a$ is a Brjuno number.

\begin{definition}
If $U\varsubsetneq\C$ is a simply connected domain containing $0$,
we denote by $\rad(U)$ the conformal radius of $U$ at $0$, i.e.,
$\rad(U)=|\phi'(0)|$ where $\phi:(\D,0)\to (U,0)$ is any conformal
representation.
\end{definition}

\begin{definition}
For any Brjuno number $\a\in {\cal B}$, we denote by $r(\a)$ the conformal
radius at 0 of the Siegel disk of the quadratic polynomial $P_\a$. If $\a\in
\R\setminus {\cal B}$, we define $r(\a)=0$.
\end{definition}

It is known that there exists a constant $C_0$ such that for any
Brjuno number $\a\in {\cal B}$ and any univalent map $f:\D\to 0$
which fixes $0$ with derivative $e^{2i\pi \a}$, $f$ has a Siegel
disk which contains $B(0,r)$ with $\Phi(\a) +\log r \geq C_0$. In
particular, for all $\a\in {\cal B}$, we have
\begin{equation}\label{yoccoz}
\Phi(\a) +\log r(\a) \geq C_0.
\end{equation}

\remark The existence of $\Delta_f$ is due to Brjuno \cite{bru}.
The lower bound (\ref{yoccoz}) is due to Yoccoz \cite{y}.
\vskip.2cm

In \cite{bc2}, we prove that there exists a universal constant
$C_1$ such that for all $\a\in {\cal B}$, we have
\begin{equation}\label{buffcheritat}
 \Phi(\a) +\log r(\a) \leq C_1.
\end{equation}

Inequalities (\ref{yoccoz}) and (\ref{buffcheritat}) imply that
$\Phi(\a)+\log r(\a)$ is uniformly bounded on ${\cal B}$:
\begin{equation}\label{ybc}
(\exists C\in \R),~(\forall \a\in {\cal B}),\quad |\Phi(\a)+\log
r(\a)|\leq C.
\end{equation}

\begin{figure}[htbp]
\centerline{
\includegraphics{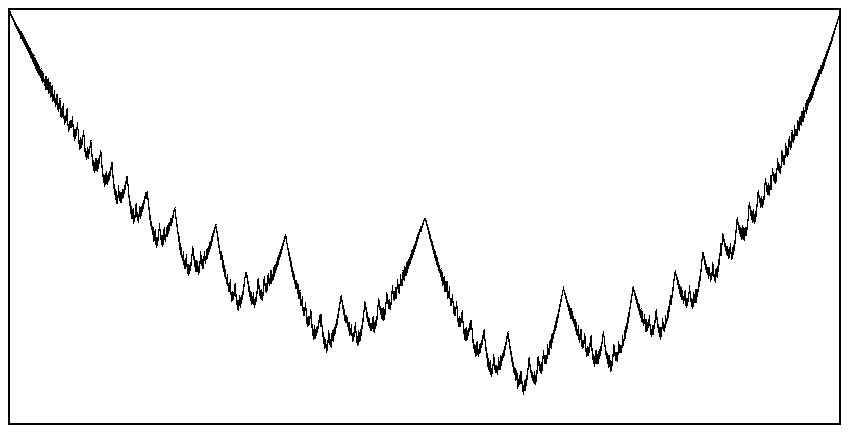}
}
\caption{The graph of the function $\a\mapsto \Phi(\a)+\log r(\a)$ with $\a\in
  [0,1]$. The range is $[0,\log(2\pi)]$. \label{graph}}
\end{figure}

In this article we prove the following result which was conjectured
by  Marmi \cite{ma}.

\begin{theorem}\label{main}
The function $\a\mapsto \Phi(\a)+\log r(\a)$ extends to $\R$ as a continuous
function.
\end{theorem}

In fact, Marmi, Moussa and Yoccoz made the following stronger
conjecture
(\cite{mmy} and \cite{ca}).

\begin{conjecture}
The function $\a\mapsto \Phi(\a)+\log r(\a)$ ---which is
well-defined on ${\cal B}$--- is H\"older of exponent $1/2$.
\end{conjecture}

\begin{figure}[htbp]
\centerline{
\includegraphics{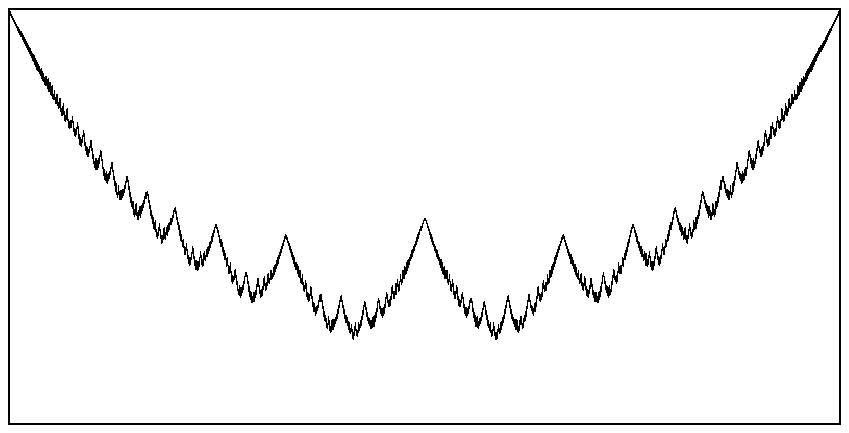}
} \caption{The graph of the function $\a\mapsto \widetilde
\Phi(\a)+\log r(\a)$
  with
  $\a\in [0,1]$. The range is
  $[0,\log(2\pi)]$. \label{graph2}}
\end{figure}

\remark In \cite{y}, Yoccoz uses a modified version of continued
fractions. He defines a sequence $\tilde\a_n$ defined by $\tilde
\a_0 = d(\a,\Z)$ and  $\tilde \a_{n+1}=d(1/\tilde \a_n,\Z)$. The
corresponding function $\widetilde \Phi$ defined by
$$\widetilde \Phi(\a)= \sum_{n\geq 0}
\tilde\a_0\cdots\tilde\a_{n-1}\log\frac{1}{\tilde\a_n}$$
has the
additional property that $\widetilde \Phi(1-\a)=\widetilde \Phi(\a)$.
Figure \ref{graph2} shows the graph of the function $\a\mapsto
\widetilde \Phi(\a)+\log r(\a)$. Theorem 4.6 in \cite{mmy} asserts
that the restriction of $\Phi-\widetilde \Phi$ to ${\cal B}$ extends
to $\R$ as a $1/2$-H\"older continuous periodic function with period
one. It follows from this result and theorem \ref{main} that the
function $\a\mapsto \widetilde \Phi(\a)+\log r(\a)$ extends to $\R$ as
a continuous function (and that the Marmi-Moussa-Yoccoz conjecture is
equivalent with $\Phi$ replaced by $\wt{\Phi})$.

\section{Statement of results.}

In this section, we will define a function $\Upsilon:\R\to \R$ and
in the rest of the article, we will show that for all $\a\in \R$,
$$\lim_{\a'\to \a,~\a'\in {\cal B}} \Phi(\a')+\log r(\a') =
\Upsilon(\a).$$
It is an easy exercise to prove that $\Upsilon$ is then continuous.

\remark For $\a\in \Q$, we give a computable formula of
$\Upsilon(\a)$. \vskip.2cm

The value of $\Upsilon$ at Brjuno numbers is obvious.

\begin{definition}\label{def1}
For $\a\in {\cal B}$, we set
$$\Upsilon(\a)=\Phi(\a)+\log r(\a).$$
\end{definition}

\subsection{Strategy of the proof.\label{sec2}}

The strategy for proving that for all $\a\in \R$,
$$\lim_{\a'\to \a,~\a'\in {\cal B}} \Phi(\a')+\log r(\a') =
\Upsilon(\a)$$
consists in bounding
$\Phi(\a')+\log r(\a')$ from above and from below as $\a'\in {\cal
B}$ tends to $\a$. The upper bound follows from techniques of
parabolic explosion developed in \cite{c} and \cite{bc2}. Those
techniques are presented in section \ref{parex}.

The lower bound essentially follows from techniques of
renormalization introduced by Yoccoz in \cite{y}. He uses
estimates which are valid for all maps which are univalent in $\D$ and
fix $0$ with derivative of modulus $1$. In our case, we will need
to improve those estimates for maps which are close to rotations
and maps which have at most one fixed point in $\D^*$ (see section
\ref{rensec}).

The rest of this section is devoted to the definition of
$\Upsilon$ at rational and \Cremer numbers.

\subsection{The value of $\Upsilon$ at rational numbers.}\label{subsec_urat}

A rational number $\a=p/q\in \Q$ has two finite continued
fraction expansions, corresponding to two sequences of approximants
$p_n/q_n$, two sequences $\alpha_n$, and two sequences $\beta_n$.
One of the sequences $\a_n$ is provided by the usual algorithm:
$\a_0 = \{\a\}$ and $\a_{n+1} = \{1/\a_n\}$, which eventually gives
$\a_m =0$ for some $m\in\N$, after which the sequence is not defined
any more. The other has the same $\alpha_k$ for $k<m$, its $\a_m=1$,
and has one more term, $\a_{m+1}= 0$.
\footnote{
A number $\a'$ tending to $p/q$ has its $\a'_k$ that
tends to the $\alpha_k$ of $p/q$ for all $k<m$. According to whether
$\a'$ tends to $p/q$ by the left or the right, $\a'_m$ tends to one
of the two values defined above, that is $0$ or $1$, the
correspondence depending on the parity of $m$. Moreover, if it is $1$,
then $\a'_{m+1}$ tends to $0$. This motivates the two
definitions we made.}

In both cases, the sequence $\beta$ is defined by $\beta_n=\a_0
\cdots \a_n$. Let $n_0 = m$ or $m+1$ be the last index of the sequence
$\a_n$ of $p/q$ that we chose. We have $\a_{n_0}=0$. We can form the
finite sum
\[\Phi_\trunc(p/q)=\sum_{n=0}^{n_0-1} \beta_{n-1}\log\frac{1}{\a_n}\]
(with the convention that a sum $\sum_{n=0}^{n=-1} \cdots$ is equal to $0$).
It turns out to be independent of the choice between the two values of
$n_0$, as can easily be checked.

%% If $\a=p/q\in \Q$ is a rational number, we can define its
%% approximants by the continued fraction algorithm and there is an
%% integer $n_0$ such that $p/q=p_{n_0}/q_{n_0}$. For $n\leq n_0$, we
%% can define $\a_n$ and $\beta_n$. Note that
%% $\a_{n_0}=\beta_{n_0}=0$.

%%\begin{definition}
%% If $\a$ is an integer, we define $\Phi_\trunc(\a)=0$.
%% If $\a=p_{n_0}/q_{n_0}\in \Q\setminus \Z$ is a rational number, we define
%% $\Phi_{\trunc}(\a)$ by
%% $$\Phi_\trunc(\a)=\sum_{n=0}^{n_0-1} \beta_{n-1}\log\frac{1}{\a_n}.$$
%%\end{definition}

The following two definitions and their relations with the
conformal radii of Siegel disks appear in \cite{c}.

\begin{definition}
Assume $f:(\C,0)\to (\C,0)$ is a germ having a multiple fixed point at the
origin whose Taylor expansion is
$$f(z) = z + Az^{q+1} + {\cal O}(z^{q+2}),\quad\text{with}\quad A\in \C^*.$$
The asymptotic size of $f$ at $0$ is defined by
$$L_a(f,0) = \left|\frac{1}{qA}\right|^{1/q}.$$
\end{definition}

\begin{definition}
Assume $p/q\in \Q$ is a rational number. Then, we define
$$L_a(p/q) = L_a(P_{p/q}^{\circ q},0).$$
\end{definition}

\begin{definition}\label{def2}
For all rational number $p/q$, we define
$$\Upsilon\left(\frac{p}{q}\right)=
\Phi_\trunc\left(\frac{p}{q}\right) +\log
L_a\left(\frac{p}{q}\right) + \frac{\log 2\pi}{q}.$$
\end{definition}

\subsection{The value of $\Upsilon$ at \Cremer numbers.}

\begin{definition}
For all irrational number $\a$ and all integer $n\geq 0$, we
define
$$\Phi_n(\a) = \sum_{k=0}^{n}\beta_{k-1}\log \frac{1}{\a_k}.$$
\end{definition}

\begin{definition}
If $U\subset\C$ is a hyperbolic connected domain containing
$0$, we denote by $\rad(U)$ the conformal radius of $U$ at $0$,
i.e., $\rad(U)=|\pi'(0)|$ where $\pi:(\D,0)\to (U,0)$ is any
universal covering.
\end{definition}

\remark This definition of conformal radius coincides with the one
given in the introduction in the case of simply connected domains.

\begin{definition}\label{def_X}
For all $\a\in \R\setminus \Q$ and all integer $n\geq 0$, we
define
$$X_n(\a) = \{z\in \C^*~|~z~\text{is a periodic point of }P_\a~\text{of
period }\leq q_n\}$$
where $p_n/q_n$ are the approximants to $\alpha$,
$$r_n(\a) = \rad(\C\setminus X_n(\a))\quad \text{and}
\quad d_n(\a) = d(0,X_n(\a)).$$
\end{definition}

\remark If $n\geq 2$, then $q_n\geq 2$, $X_n(\a)$ contains at least two
points and $r_n(\a)\in ]0,\infty[$. Moreover, for $n\geq 2$, the
function $\a\mapsto \log r_n(\a)$ is well-defined and continuous
in a neighborhood of every point $\a\in \R\setminus \Q$.

\vskip.2cm

For all irrational number $\a$, the sequence $(r_n(\a))_{n\geq 0}$
is decreasing and converges to $r(\a)$ as $n\to\infty$. Indeed, if $0$
is not linearizable, it is accumulated by periodic points of
$P_\a$.\footnote{In fact, Yoccoz
proved that $0$ is accumulated by whole cycles.}
If $0$ is linearizable, the Siegel disk $\Delta_\a$ is contained
in $\C\setminus X_n(\a)$ for all $n\geq 0$ and the boundary of
$\Delta_\a$ is accumulated by periodic points of $P_\a$.\footnote{It
is not known whether $\partial \Delta_\alpha$ is always accumulated by
whole cycles.}
Since $P_\alpha$ is tangent the rotation of angle $\alpha$
and $\a$ is irrational, if $0$ is not linearizable,
then \[r_n(\alpha) \underset{n\to +\infty}{\sim} d_n(\alpha).\]
If $\a$ is a Brjuno number, then
$$\lim_{n\to \infty} \Phi_n(\a)+\log r_n(\a)=\Upsilon(\a).$$

Yoccoz's work \cite{y} implies that there exists a constant $C_0$ such
that for all $\a\in \R\setminus \Q$ and all $n\geq 0$,
$$\Phi_n(\a)+\log r_n(\a)\geq C_0.$$
Thus, the following definition makes sense.

\begin{definition}\label{def3}
For all \Cremer number $\a$, we define
$$\Upsilon(\a) = \liminf_{n\to \infty} \Phi_n(\a)+\log r_n(\a).$$
\end{definition}

In fact, we will see (section \ref{parex}) that this is an actual
limit.

\begin{theorem}\label{limitcremer}
For all \Cremer number $\a$,
$$\Upsilon(\a)= \lim_{n\to +\infty} \Phi_n(\a)+\log r_n(\a).$$
\end{theorem}

\begin{corollary}
For all \Cremer number $\a$,
$$\Upsilon(\a)= \lim_{n\to +\infty} \Phi_n(\a)+\log d_n(\a).$$
\end{corollary}

Our goal is to prove that for all $\a\in \R$, the value of
$\Upsilon(\a)$ defined previously (see definitions \ref{def1},
\ref{def2} and \ref{def3}) is the limit of $\Phi(\a')+\log r(\a')$
as $\a'\in {\cal B}$ tends to $\a$. In section \ref{parex} we
introduce the techniques of parabolic explosion and in section
\ref{prooflimsup} we show that for all $\a\in \R$,
\begin{equation}\label{limsup}
\limsup_{\a'\to \a,~\a'\in {\cal B}} \Phi(\a')+\log r(\a') \leq
\Upsilon(\a).
\end{equation}
In section \ref{rensec}, we refine Yoccoz's estimates for
renormalization of univalent maps $f:\D\to \C$ which fix $0$ with
derivative of modulus $1$, and in sections \ref{proofliminf1} and
\ref{cremerparabo} we show that for all $\a\in \R$,
\begin{equation}\label{liminf}
\liminf_{\a'\to \a,~\a'\in {\cal B}} \Phi(\a')+\log r(\a') \geq
\Upsilon(\a).
\end{equation}
Let us mention that inequality (\ref{limsup}) without inequality
(\ref{liminf}) (respectively inequality (\ref{liminf}) without
inequality (\ref{limsup})) is not sufficient to conclude that
$\Upsilon$ is upper semi-continuous (respectively lower
semi-continuous) since we only consider approximating $\a$ with
sequences of Brjuno numbers.

\section{Parabolic explosion.\label{parex}}

In this section, we first present the techniques of parabolic
explosion. We then apply those techniques in order to prove
theorem \ref{limitcremer}.

\subsection{Definitions.}

Assume $p/q\in \Q$ is a rational number. The origin is a parabolic fixed point
for the quadratic polynomial $P_{p/q}$. It is  known (see \cite{dh},
chapter IX) that there exists a complex number $A\in \C^*$ such that
\[P_{p/q}^{\circ q}(z) = z+Az^{q+1} + {\cal O}(z^{q+2}).\]
Thus, $P_{p/q}^{\circ q}$ has a fixed point of multiplicity $q+1$ at the
origin. By Rouché's theorem, when $\a$ is close to $p/q$, the polynomial
$P_\a^{\circ q}$ has $q+1$ fixed points close to $0$. One coincides with
$0$. The others form a cycle of period $q$ for $P_\a$. More precisely, we have
the following (see \cite{c} or \cite{bc2} proposition 1 for a proof).

\begin{proposition}\label{functionchi}
Let $p/q$ be a rational number, and $\zeta= e^{2i\pi p/q}$.
There exists an analytic function
$\chi:B(0,1/q^{3/q})\to \C$ such that
$\chi(0)=0$ and for any $\delta\in
B(0,1/q^{3/q})\setminus\{0\}$, $\chi(\delta)\neq 0$ and the set
$$\Big< \chi(\delta),\chi(\zeta\delta),
\chi(\zeta^2\delta),\ldots,\chi(\zeta^{q-1}\delta)\Big>$$
forms a cycle of period $q$ of $P_{p/q+\delta^q}$.
We will note $\chi = \chi_{p/q}$, since it depends on $p/q$.
\end{proposition}

\remark Observe that $\d\in B(0,1/q^{3/q})$ if and only if
$\a=p/q+\d^q\in B(p/q,1/q^3)$. \vskip.2cm

In the following definition, note that $\a$ is a {\bf complex} number.

\begin{definition}
For all $p/q\in \Q$ and all $\a\in B(p/q,1/q^3)$, we define
$${\cal C}_{p/q}(\a)=\chi_{p/q}\left\{\sqrt[q]{\a-p/q}\right\},$$
where $\sqrt[q]{z}$ denotes the \emph{set} of complex $q$-th roots
of $z$.
\end{definition}

The set ${\cal C}_{p/q}(\a)$ is a cycle of period $q$ for $P_\a$,
except when $\a = p/q$, in which case it is reduced to $\{0\}$.
In particular, if $\a$ is irrational, $p/q=p_n/q_n$ is an
approximant to $\a$ and $|\a-p_n/q_n|<1/q_n^3$, then
${\cal C}_{p_n/q_n}(\a)\subset X_n(\a)$. 
Note that when $|\a_0-p/q|<1/2q^3$, the cycle ${\cal C}_{p/q}(\a)$
is defined for all $\a\in B(\a_0,1/2q^3)$, and not reduced to $\{0\}$.

\subsection{A preliminary lemma.}

\begin{lemma}\label{multiple}
Assume $\a_0\in \R\setminus \Q$ and let $p_n/q_n$ be an
approximant to $\a_0$ with $q_n\geq 2$. Assume $\alpha\in\C$, $\a
\neq p_n/q_n$, $q \leq q_n$ and $P_\alpha^{\circ q}$ has a
multiple fixed point. Then, $$|\alpha_0-\alpha| \geq
\frac{1}{2q_n^3}.$$
\end{lemma}

\proof Either $\a=p/q$ for some integer $p$. Within the disk
$B(\a_0,1/2q_n^3)$, the only possibility is $p/q=p_n/q_n$. Or $\a$
belongs to a Yoccoz disk of radius $\log 2/(2\pi q') < 1/8q'$
tangent to the real axis at $p'/q'$ for some rational number
$p'/q'$ with $q'<q\leq q_{n}$. By a well-known property of
approximants, we have
$$|q'\a_0-p'|\geq|q_{n-1}\a_0+p_{n-1}| \geq \frac{1}{q_n+q_{n-1}}
\geq \frac{1}{2q_{n}}.$$ Moreover, by Pythagoras' theorem,
\begin{eqnarray*}
|\a-\a_0| & \geq & \frac{1}{q'}\left(\sqrt{\left(q'\a_0-p'\right)^2+
 \left(1/8\right)^2}-1/8\right)\\
& \geq & \frac{1}{q_n}
 \left(\sqrt{1/\left(2q_n\right)^2+
 1/8^2}-1/8\right) \\
& = & \frac{1/(2q_n)^2}{q_n \Big(
 \sqrt{1/\left(2q_n\right)^2+1/8^2}+1/8^2\Big)}\\
& \geq & \frac{1}{2q_n^3}\cdot \frac{1}{2 \Big( \sqrt{1/4^2 + 1/8^2} +
 1/8^2\big)} \quad \geq \frac{1}{2q_n^3}.
\end{eqnarray*}
 \qedprop

\begin{corollary}\label{holoq3}
Assume $\a_0\in \R\setminus \Q$ and let $p_n/q_n$ be an
approximant to $\a_0$ with $q_n\geq 2$. The set
$$X(\a)=\{z\in \C^*~|~z~\text{is a periodic point of}~P_\a~\text{of period
}\leq q_n\}$$ moves holomorphically with respect to $\a\in
B(\a_0,1/2q_{n+1}^3)$.
\end{corollary}

\proof If the set $X(\a)$ fails to move holomorphically at a point
$\a\in\C$, then, for some integer $q\leq q_{n}$, $P_\a^{\circ q}$
has a multiple fixed point. Either $\a = p_n/q_n$, and (according
to a property of approximants) $|\a-\a_0| \geq 1/(2q_n q_{n+1}) >
1/2q_{n+1}^3$. Or $\a \neq p_n/q_n$, and by the previous lemma
$|\a-\a_0| \geq 1/2q_n^3 > 1/2q_{n+1}^3$. \qedprop

\subsection{A technical lemma.}

\begin{lemma}\label{lemmatec}
There exists $C\in \R$ such that for all $\a_0\in \R\setminus \Q$
and all $p/q\in \Q$ with $q\geq 2$, the following holds. Assume
$V(\a)\ni 0$ is an open set that moves holomorphically with
respect to $\a\in B(\a_0,1/2q^3)$.
\begin{itemize}
\item If $|\a_0-p/q|\geq 1/2q^3$, set $V'(\a_0)=V(\a_0)$.
\item If
$|\a_0-p/q|< 1/2q^3$, assume ${\cal C}_{p/q}(\a)\subset V(\a)$ for
all $\a\in B(\a_0,1/2q^3)$ and set $V'(\a_0)=V(\a_0)\setminus
{\cal C}_{p/q}(\a_0)$.
\end{itemize}
Then,
$$\log \frac{\rad(V'(\a_0))}{\rad(V(\a_0))} \leq \frac{\log |\a_0-p/q|}{q} +
C\frac{\log q}{q}.$$
\end{lemma}

\proof Let us first assume that $|\a_0-p/q|\geq 1/2q^4\geq 1/q^5$
(this comprises the case $V'(\alpha_0) = V(\alpha_0)$).
Then,
$$\log |\a_0-p/q| + 5\log q\geq 0$$ and the
lemma follows trivially with $C=5$ since
$$\log \frac{\rad(V'(\a_0))}{\rad(V(\a_0))} \leq 0.$$ So, let us assume that
$|\a_0-p/q|<1/2q^4.$ Then,
\[B\egaldef B(p/q,1/2q^4)\subset
B(\a_0,1/q^4)\subset B(\a_0,1/2q^3).\] We set
$$U=\{\d\in \C~|~p/q+\d^q\in B \}
\quad\text{and}\quad
S=\{\delta\in U~|~p/q+\delta^{q}=\a_0\}.$$
Note that $\chi_{p/q}(S)={\cal C}_{p/q}(\a_0)$.

The radius of the disk $U$ is $1/(2q^4)^{1/q}$ and the set $S$
consists in $q$ points equidistributed on a circle of radius
$|\a_0-p/q|^{1/{q}}.$ So, according to
proposition~\ref{prop_BC2_12} (see the appendix~\ref{app_A}), we
have
$$\log\frac{\rad(U\setminus S)}{\rad(U)}  <
\log \frac{|\a_0-p/q|^{1/q}}{1/(2q^4)^{1/q}} + \frac{C}{q}$$ for some
universal constant $C$.

According to proposition~\ref{prop_BC2_13} (see the
appendix~\ref{app_A}), there exists for $\a\in B(\a_0,1/2q^3)$ an
analytic family of universal coverings $\pi_\a : \widetilde{V}(\a)
\to V(\a)$, where $\widetilde{V}(\a)$ are open subsets of
$B(0,4)$, and $\widetilde{V}({\a_0}) = \D$. The set $V(\a)$ moves
holomorphically with $\a\in B(\a_0,1/2q^3)$ and when $\d\in U$,
$\a(\d)=p/q+\d^q$ belongs to $B\subset B(\a_0,1/q^4)$. For $\a \in
B$, the sets $\widetilde{V}(\a)$ are all contained in some ball
$B(0,\rho)$ with
$$\log \rho= \frac{2\log 4}{1+\ds\frac{1/2q^3}{1/q^4}} = \frac{\log 16}{1+q/2}.$$

The map $\chi_{p/q}$ ``lifts'' to a map
$\phi : U \to B(0,\rho)$ such that $\phi(\delta) \in \widetilde{V}(\a(\delta))$.
It follows from the definitions that,
$$\log \frac{\rad(V'(\a_0))}{\rad(V(\a_0))}=\log
  \frac{\rad(V(\a_0)\setminus {\cal 
  C}_{p/q}(\a_0))}{\rad(V(\a_0))}
= \log\frac{\ds \rad\left(\widetilde{V}(\a_0) \setminus
  \pi_{\a_0}^{-1}\left(\chi_{p/q}(S)\right)\right)}
{\rad(\widetilde{V}(\a_0))}.$$
Now $\widetilde{V}(\a_0) = \D$ and $\phi(S) \subset
\pi_{\a_0}^{-1}(\chi_{p/q}(S))$, thus
$$\log \frac{\rad(V'(\a_0))}{\rad(V(\a_0))}
\leq \log \rad(\D \setminus \phi(S))\leq \log \rad(B(0,\rho)
\setminus \phi(S)).$$ The range of the function $\phi$ needs not to be
a subset of $\D$, but we know 
proposition~\ref{prop_BC2_10} (see the
appendix~\ref{app_A}), 
\begin{eqnarray*}
\log \rad(B(0,\rho) \setminus \phi(S)) & \leq &
\log\frac{\rad(U\setminus S)}{\rad(U)}+\log \rho\\
& \leq & \frac{\log|\a_0-p/q|}{q}+4\frac{\log
  q}{q}+\frac{\log 2}{q}+\frac{C}{q}+\frac{\log 16}{1+q/2}\\
& \leq & \frac{\log|\a_0-p/q|}{q}+C'\frac{\log q}{q}
\end{eqnarray*}
for some universal constant $C'$. \qedprop

\subsection{A short remark}

Let $F_n$ be the smallest possible value of
$q_n$ over all irrationals $\a$, where $p_n/q_n$ is the $n$-th
approximant to $\a$. Then 
$F_n$ is the Fibonacci sequence defined by
\[F_{-1} = 0,\ F_0 = 1,\ F_{n+1} = F_{n}+F_{n-1}.\]
The first terms are
\[F_{-1} = 0,\ F_0 = 1,\ F_1=1,\ F_2=2,\ F_3=3,\ F_4=5,\ \ldots\]
The function $x \mapsto x/\log x$ is decreasing on $[e,+\infty[$, thus
\[\forall n\geq 3,\ \frac{\log q_n}{q_n} \leq \frac{\log F_n}{F_n}.\]
For $n=1$ and $2$, the biggest possible value of $\log
(q_n)/q_n$ is $\log(3)/3$.

\subsection{An important corollary.}

The next proposition tells us that for all irrational $\a$, the
sequence $\Phi_n(\a)+\log r_n(\a)$ is essentially decreasing, in the
sense that it can not increase too fast.

\begin{proposition}\label{decprop}
There exists a constant $C\in \R$ such that for all
$\a\in\R\setminus \Q$ and all $n\geq 1$ such that $q_{n}\geq 2$
(with $p_n/q_n$ the approximants to $\alpha$), we have
$$\Big(\Phi_{n+1}(\a)+\log r_{n+1}(\a)\Big)-\Big(\Phi_{n}(\a)+\log
r_{n}(\a)\Big)
\leq C\frac{\log q_{n+1}}{q_{n+1}}.$$
% \leq C\frac{\log F_{n+1}}{F_{n+1}} %% On peut pas à cause du cas n=1...
\end{proposition}

\proof Let us fix $\a_0\in \R\setminus \Q$ and choose $n$ so that
$q_n\geq 2$. We want to apply lemma \ref{lemmatec} with
$p/q=p_{n+1}/q_{n+1}$ and
$$V(\a)= \C\setminus \{z\in \C^*~|~z~{\rm
is~a~periodic~point~of~}P_\a~\text{of period }\leq q_n\}.$$ By
definition, $0\in V(\a)$ and by corollary \ref{holoq3}, the set
$V(\a)$ moves holomorphically with respect to $\a\in
B(\a_0,1/2q_{n+1}^3)$. Also, $V(\a)$ contains the periodic cycles
of $P_\a$ of period $q_{n+1}$ and so, if $|\a_0-p/q|<1/2q^3$, then
${\cal C}_{p/q}(\a)\subset V(\a)$ for all $\a\in B(\a_0,1/2q^3)$.
As in lemma \ref{lemmatec}, if $|\a_0-p/q|\geq 1/2q^3$, we set
$V'(\a_0)=V(\a_0)$  and otherwise, we set
$V'(\a_0)=V(\a_0)\setminus {\cal C}_{p/q}(\a_0)$. Then,
$$r_n(\a_0) = \rad(V(\a_0))\quad\text{and}\quad
r_{n+1}(\a_0) \leq \rad(V'(\a_0)).$$ So, lemma \ref{lemmatec}
implies that
\begin{eqnarray*}
\log r_{n+1}(\a_0)- \log r_n(\a_0) & \leq & \frac{\log
|\a_0-p_{n+1}/q_{n+1}|}{q_{n+1}}+C\frac{\log q_{n+1}}{q_{n+1}}
\\
& = & \frac{\log \beta_{n+1}}{q_{n+1}} +(C-1)\frac{\log
  q_{n+1}}{q_{n+1}}.
\end{eqnarray*}
Since $\beta_{n+1} \leq \alpha_{n+1}$ and $1/q_{n+1} \geq \beta_n$:
\begin{eqnarray*}
\log r_{n+1}(\a_0)- \log r_n(\a_0) & \leq
& -\beta_n\log \frac{1}{\a_{n+1}} +(C-1)\frac{\log q_{n+1}}{q_{n+1}}\\
& = & -\Phi_{n+1}(\a_0)+\Phi_n(\a_0) + (C-1)\frac{\log
q_{n+1}}{q_{n+1}}
\end{eqnarray*}
for some universal constant $C$. \qedprop

The bound we gave depends on $\a$, but for each $n$, the supremum
over all $\a \in \R\setminus\Q$ is exponentially decreasing.

\subsection{Application to the proof of theorem \ref{limitcremer}.}

Assume $\a$ is a \Cremer number, define $u_n=\Phi_n(\a)+\log
r_n(\a)$ and let us recall that by definition, $\Upsilon(\a)=
\liminf_{n\to \infty} u_n.$ The sequence $u_n$ is not decreasing, but
it is ``essentially decreasing'', in the sense that
proposition \ref{decprop} gives us
$$u_{n+1}-u_n \leq C\frac{\log q_{n+1}}{q_{n+1}}$$
and $(\log q_{n+1})/q_{n+1}$ decreases
exponentially fast.
Therefore the sequence $u_n$ converges: indeed, if we choose $n_0$
large enough so that 
$$\sum_{n\geq n_0} C\frac{\log q_{n+1}}{q_{n+1}}\leq \e
\quad\text{and}\quad u_{n_0}\leq \Upsilon(\a) +\e,$$ then
$\Upsilon(\a) \leq u_n\leq \Upsilon(\a) +2\e$ for all $n\geq n_0$.

\section{Proof of inequality (\ref{limsup}).\label{prooflimsup}}

\subsection{Irrational numbers.}

We will now show that for all $\a\in \R\setminus \Q$,
$$\limsup_{\a'\to \a,~\a'\in {\cal B}} \Phi(\a')+\log r(\a') \leq
\Upsilon(\a).$$ Let us fix $\e>0$. We must show that for $\a'\in
{\cal B}$ sufficiently close to $\a$, $\Phi(\a')+\log r(\a')\leq
\Upsilon(\a)+\e$. Remember that as $n\to \infty$, $\Phi_n(\a)+\log
r_n(\a)\to \Upsilon(\a)$. So, let us choose $n_0$ large enough so
that
$$\Phi_{n_0}(\a)+\log r_{n_0}(\a)\leq \Upsilon(\a)+\e/3.$$
Increasing $n_0$ if necessary, we may also assume that 
%for all $\a'\in \R\setminus \Q$ with approximants $p'_n/q'_n$, we have
%$$\sum_{n\geq n_0} C\frac{\log q'_{n+1}}{q'_{n+1}}\leq \e/3,$$
$n_0 \geq 2$ and
$$\sum_{n\geq n_0} C\frac{\log F_{n+1}}{F_{n+1}}\leq \e/3,$$
where $C$ is the constant in proposition
\ref{decprop}. In a neighborhood of $\a$, the
functions $\Phi_{n_0}$ and $\log r_{n_0}$ are continuous. So, if
$\a'$ is sufficiently close to $\a$,
$$\Phi_{n_0}(\a')+\log r_{n_0}(\a')\leq \Phi_{n_0}(\a)+\log
r_{n_0}(\a)+\e/3$$ and summing the inequality of proposition
\ref{decprop} from $n=n_0$ to $n=+\infty$ yields
$$\Phi(\a')+\log
r(\a')\leq \Upsilon(\a)+\e.$$
\qedprop

\subsection{Rational numbers.\label{parexparabo}}

We will
show that $$\limsup_{\a'\to p/q,~\a'\in {\cal B}} \Phi(\a')+\log
r(\a') \leq \Upsilon(p/q).$$
In the whole section, we will use the notation $$\e=\a'-p/q.$$
For $\a'\in\C$ and $\theta\in\R$, we will also denote
$R_{\a'}(\theta)$ the external ray of argument $\theta$ of
$P_{\a'}$. The external rays for the Mandelbrot set will be noted
$R_M(\theta)$.

The polynomial $P_{\a}$ is conjugate to the quadratic polynomial
$z\mapsto z^2+c$ with $c=e^{2i\pi \a}/2- e^{4i\pi\a}/4$.
When $\Im(\a) \tend -\infty$ and $\Re(\a) \tend \wt{\theta}$,
then $|c|\tend +\infty$ and $\on{arg} c \tend 2 \wt{\theta} + \frac{1}{2}
\bmod 1$. Given $\wt{\theta}\in\R$, we will denote by ${\cal R}(\wt{\theta})$
the connected component of 
the preimage of $R_M(2\wt{\theta}+1/2)$ by
$\a \mapsto c$, whose real part tends to $\wt{\theta}$.

When $\a$ is real,
the parameter $c$ is on the boundary of the main cardioid of the
Mandelbrot set. If $\a = p/q\notin \Z$, $c\neq 1/4$ and there are two
external rays of $M$ landing at $c$. We denote by
$\theta^-<\theta^+$ their arguments in $]0,1[$. The arguments
$\theta^+$ and $\theta^-$ are periodic of period $q$ under
multiplication by $2$ modulo $1$. They belong to the same orbit
$\Theta$. In the dynamical plane of $P_{p/q}$, the rays
$R_{p/q}(\theta)$, $\theta\in \Theta$, form a periodic cycle of
rays which land at $0$. If $p/q\in \Z$, the dynamical ray of
argument $0$ is fixed and lands at $0$. We set
$\theta^-=\theta^+=0$ and $\Theta=\{0\}$.

Let us recall the following rule: the ray
$R_{\a'}(\theta)$ moves holomorphically with $\a'$ as long as
$c$ does not belong to the closure of the union of the $R_M(2^k
\theta)$ for $k\in\N^*$.

\begin{definition}
When $\a'\in \R$ is close to $p/q$,
the rays $R_{\a'}(\theta)$, $\theta\in \Theta$, form a
cycle of rays which land on the cycle ${\cal
C}_{p/q}(\a')$. We denote by $Y(\a')$ the union of ${\cal
C}_{p/q}(\a')$ and this cycle of rays.
\end{definition}

\begin{figure}[htbp]
\centerline{ \fbox{\scalebox{.9}{\includegraphics{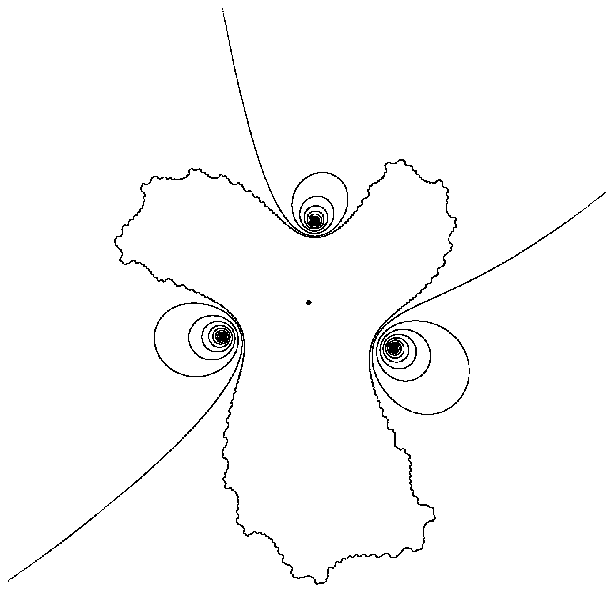}}}
\quad \fbox{\scalebox{.9}{\includegraphics{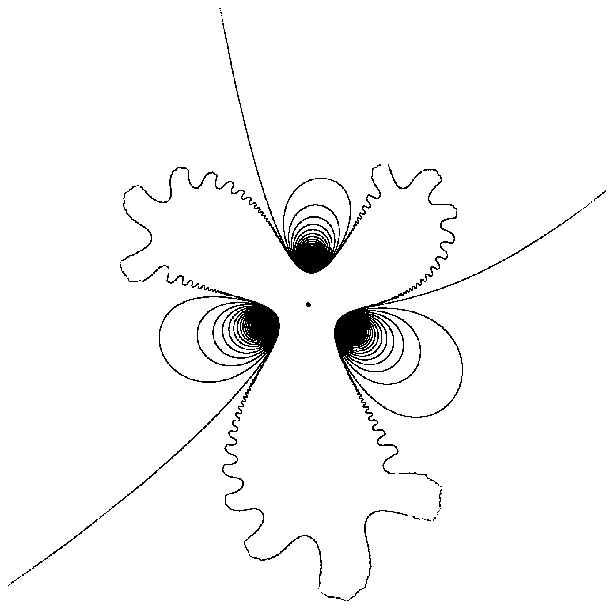}}} }
\caption{The rays of argument $1/7$, $2/7$ and $4/7$ and the
boundary of the Siegel disk for the polynomial $P_{(1/3)+\e}$: left
for $\e=\sqrt 2/1000$ and right for $\e=\sqrt
2/10000$.\label{rays}}
\end{figure}

Figure \ref{rays} shows the rays of argument $1/7$, $2/7$ and
$4/7$ and the boundary of the Siegel disk for the polynomial
$P_{(1/3)+\e}$ for $\e=\sqrt 2/1000$ and $\e=\sqrt 2/10000$.

If $\e$ is irrational and is close enough to $0$, then
$p/q$ is an approximant $p'_{n_0}/q'_{n_0}$ to $\a'$,
and its index $n_0$
is the same number as in section~\ref{subsec_urat} and depends on
the sign of $\e$.
% Note that $n_0$ is either equal to $m$ or $m+1$,
% depending on the sign of $\e$.
As $\a'\to p/q$, $\log 
\rad(\C\setminus Y(\a'))\to -\infty$ and $\beta'_{n_0-1}\log
(1/\a'_{n_0})\to -\infty$.
% In fact, those two quantities compensate each other.
We postpone the proof of the following lemma to section \ref{secya}.
% The following lemma is proved in section \ref{secya}.

\begin{lemma}\label{sizeya}
We have
$$\limsup_{\a'\to p/q,~\a'\in \R\setminus \Q}
\log \rad(\C\setminus Y(\a'))+
\beta'_{n_0-1}\log\frac{1}{\a'_{n_0}} \leq  \log
L_a\left(\frac{p}{q}\right)+\frac{\log 2\pi}{q}.$$
\end{lemma}

When $\a'$ is close to $p/q$ but not necessarily real, the
dynamical rays of argument $\theta\in \Theta$ may bifurcate. In a
neighborhood of $p/q$, this precisely occurs when $c'=e^{2i\pi
\a'}/2- e^{4i\pi \a'}/4$ belongs to ${\cal R}_M(\theta^+)$ or
${\cal R}_M(\theta^-)$.

\begin{lemma}\label{crat}
There exists a constant $c\in ]0,1]$, which depends on $p/q$,
such that the following holds. Assume $\a'\in \R\setminus \Q$ and
$p/q$ is an approximant to $\a'$. Let $n_0$ be its index. Let
$p'_{n_0+1}/q'_{n_0+1}$ be $\a'$'s next approximant. 
Then, for all $\a''\in
B(\a',c/(q'_{n_0+1})^2)$, the dynamical rays of argument
$\theta\in \Theta$ do not bifurcate. In particular, $Y(\a'')$
moves holomorphically with respect to $\a''\in
B(\a',c/(q'_{n_0+1})^2)$.
\end{lemma}

\proof
There is exactly one pair $\wt{\theta}^- < \wt{\theta}^+$,
with $2\wt{\theta}^++1/2 = \theta^+$ and $2\wt{\theta}^-+1/2 = \theta^-$
such that
${\cal R}(\wt{\theta}^+)$ and ${\cal R}(\wt{\theta}^-)$ land on $p/q$. 
The rays ${\cal R}(\wt{\theta}^+)$ and ${\cal R}(\wt{\theta}^-)$ are separated
from the upper half plane (that corresponds to 
the cardioid by $\a \mapsto c$), by a smooth curve having a
contact of order 2 with the real line, at $p/q$.
Also, the other external rays $R_M(\theta')$ for $\theta' \in
\Theta\setminus\{\theta^+, \theta^-\}$ do not land on the cardioid.
Therefore, there exists a
constant $c'>0$ such that the dynamical rays of argument
$\theta\in \Theta$ do not bifurcate when $\a''\in
B(\a',c'|\a'-p/q|^2)$. The result follows since
$$\left|\a'-\frac{p}{q}\right|^2 \geq
\left(\frac{1}{2q'_{n_0}q'_{n_0+1}}\right)^2=
\frac{1}{4q^2 (q'_{n_0+1})^2}. $$ \qedlem

Let us choose $c$ as in lemma \ref{crat} and $\a'\in {\cal B}$
sufficiently close to $p/q$ so that
$q'_{n_0+1}>1/2c$ (we denote by $p'_n/q'_n$ the approximants to
$\a'$). Then, the set $Y(\a'')$ moves holomorphically with respect
to $\a''\in B(\a',1/2(q'_{n_0+1})^3)$. Let us also assume that
$q'_{n_0+1}\geq 2$

\begin{lemma}\label{equarat}
Under the assumptions above, we have
$$\Phi(\a')+\log r(\a')\leq \Phi_{n_0}(\a')+\log
\rad(\C\setminus Y(\a')) + (C-1)\sum_{n\geq n_0+1}\frac{\log
q'_n}{q'_n},$$ where $C$ is the constant provided by lemma
\ref{lemmatec}. \end{lemma}

\proof For $\a''\in B(\a',1/2(q'_{n_0+1})^3)$, let us define
$V_{n_0}(\a'')=\C\setminus Y(\a'')$ and by induction, for $n\geq
n_0+1$ and $\a''\in B(\a',1/2(q'_{n+1})^3)$, let us define
\begin{itemize}
\item $V_{n}(\a'') = V_{n-1}(\a'')\setminus {\cal
C}_{p'_n/q'_n}(\a'')$ if $|\a'-p'_n/q'_n|<1/2(q'_n)^3$ and \item
$V_{n}(\a'') = V_{n-1}(\a'')$ otherwise.
\end{itemize}
Then, the
hypotheses of lemma~\ref{lemmatec} are satisfied and (as in
proposition \ref{decprop}), we have
\begin{eqnarray*}
\log \rad(V_{n}(\a')) - \log \rad(V_{n-1}(\a')) & \leq &
\frac{\log |\a'-p'_n/q'_n|}{q'_n} + C\frac{\log q'_n}{q'_n}\\
& \leq & -\Phi_n(\a')+\Phi_{n-1}(\a') + (C-1)\frac{\log
q'_n}{q'_n},
\end{eqnarray*}
where $C$ is the constant provided by lemma \ref{lemmatec}. The
Siegel disk $\Delta_{\a'}$ is contained in the intersection of the
sets $V_n(\a')$, and so,
$$\log r(\a') - \log \rad(V_{n_0}(\a')) \leq
-\Phi(\a') + \Phi_{n_0}(\a') + (C-1)\sum_{n\geq n_0+1}\frac{\log
q'_n}{q'_n}.$$ \qedlem

As $\a'$ tends to $p/q$, each $q'_{n_0+k}$ (for $k\geq 1$) tends to $\infty$,
thus the $n_0+k$-th summand tends to $0$. Since the sum is dominated by a
summable sequence ($\log (F_n)/F_n$), this yields
$$\sum_{n\geq n_0+1}\frac{\log q'_n}{q'_n}\to 0.$$
%(the sequence$\log q'_n/q'_n$ decreases universally exponentially fast).
Moreover, $\Phi_{n_0-1}(\a')$ converges to $\Phi_\trunc(p/q)$ and
by lemma \ref{sizeya},
$$\limsup_{\a'\to p/q,~\a'\in \R\setminus \Q}
\Phi_{n_0}(\a') + \log \rad(\C\setminus Y(\a'))
\leq \Upsilon(p/q).$$ This completes the proof of inequality
(\ref{limsup}).

\subsection{Proof of lemma \ref{sizeya}.\label{secya}}

%be the following periodic point of $P_{\a}$:
%$$\left\{\begin{array}{rl}
%z_\e & = \chi_{p/q}(\sqrt[q]\e)\quad\text{if}~\e>0,\\
 %z_\e & = \chi_{p/q}(e^{i\pi/q}\sqrt[q]{-\e})\quad\text{if}~\e<0.
%\end{array}\right.$$

We recall that $\a' = p/q+\epsilon$ is real, and that $n_0$
depends on the sign of $\epsilon$. 

\begin{lemma}\label{moduleze}
For $\e\in \R^*$ small enough, let $z_\e$ be a periodic point of
$P_{\a'}$ in the cycle ${\cal C}_{p/q}(\a')$. Then,
$$\log|z_\e| + \beta'_{n_0-1}\log\frac{1}{\a'_{n_0}}= \log
L_a\left(\frac{p}{q}\right) + \frac{\log 2\pi}{q} + {\cal
O}(\e^{1/q}).$$
\end{lemma}

\proof
By definition of the asymptotic size, we have
$$L_a(p/q) = \left|\frac{1}{qA}\right|^{1/q}
\quad\text{with}\quad
P_{p/q}^{\circ q}(z) = z+Az^{q+1}+{\cal O}(z^{q+2}).$$
Moreover, $P_{p/q+\e}^{\circ q}(0)=0$ and $(P_{p/q+\e}^{\circ
  q})'(0)=e^{2i\pi q\e}$. So
$$P_{p/q+\e}^{\circ q}(z) = e^{2i\pi q\e}z+Az^{q+1}+{\cal O}(\e z^2).$$
We know that $z_\epsilon \tend 0$ and that $P_{p/q+\e}^{\circ q}(z_\e) = z_\e$.
Therefore, we have
$$z_\e^q = \frac{1-e^{2i\pi q\e}}{A}\left( 1+ {\cal O}(z_\e)\right)
=\frac{-2i\pi q\e}{A}\left( 1+ {\cal O}(z_\e) + {\cal O}(\e)\right).$$
Thus, $z_\e={\cal O}(\e^{1/q})$ and
$$\log|z_\e| = \frac{1}{q}\log\left|\frac{2\pi q \e}{A}\right|+{\cal
  O}(\e^{1/q}).$$
Observe that
$$\frac{1}{q}\log\left|\frac{2\pi q \e}{A}\right|=
\log L_a\left(\frac{p}{q}\right) + \frac{\log 2\pi}{q} +
\frac{1}{q}\log q^2|\e|.$$
Now, if $\a'$ is sufficiently close to
$p/q$, then the $n_0$-th approximant $p'_{n_0}/q'_{n_0}$ to $\a'$
is $p/q$, and therefore when $\e$'s sign is fixed, $n_0$ is fixed,
and the numbers $q'_{n_0}$ and $q'_{n_0-1}$ are constants.
We have
\begin{eqnarray*}
\beta'_{n_0-1} & = & |q'_{n_0-1}\a' -p'_{n_0-1}| = \left|q'_{n_0-1}
\Big(\frac{p_{n_0}}{q_{n_0}} + \e\Big) - p'_{n_0-1}\right|\\
 & = & \left|\frac{1}{q'_{n_0}}\pm q'_{n_0-1}\e\right| =\frac{1}{q'_{n_0}} +
{\cal O}(\e), \quad \text{and}\\
\beta'_{n_0} & = & q'_{n_0}|\e|, \quad \text{thus}\\
\alpha'_{n_0} & = & \frac{\beta'_{n_0}}{\beta'_{n_0-1}} = (q'_{n_0})^2|\e|(1
+{\cal O}(\e)).
\end{eqnarray*}
%% Thus, we have
%% $$\frac{1}{q}\log q^2|\e| = \frac{1}{q'_{n_0}} \log |\a'_{n_0}| + {\cal O}(\e)
%% = \beta'_{n_0-1}\log |\a'_{n_0}| +{\cal O}(\e\log|\e|).$$ \qedprop
Thus, we have
$$\beta'_{n_0-1}\log |\a'_{n_0}| = \Big(\frac{1}{q} + {\cal
O}(\e)\Big)\log\big(q^2 |\e| (1+ {\cal O}(\e))\big)
= \frac{1}{q}\log q^2|\e| +{\cal O}(\e\log|\e|).$$ \qedprop

Let us now study the dynamical behaviour of $P_{p/q+\e}$ at the
scale of $z_\e$. For this purpose, we rescale the dynamical plane.
More precisely, we introduce the conjugate polynomial
$$Q_\e:w\mapsto \frac{1}{z_\e}P_{p/q+\e}(z_\e w).$$
This polynomial is conjugate to $P_{p/q+\e}$. It fixes $0$ with derivative
$e^{2i\pi (p/q+\e)}$ and has a cycle of period $q$ containing $1$.

As $\e\to 0$, $Q_\e$ converges uniformly on every compact subset of $\C$ to
the rotation $w\mapsto e^{2i\pi p/q}w$. Hence, $Q_\e^{\circ q}$ converges
uniformly on every compact subset of $\C$ to the identity. However, the
limit of the dynamics of $Q_\e$ is richer than the dynamics of the identity.
In some sense, it contains the real flow of the vector field
$2i\pi qw(1-w^q)\frac{\partial}{\partial w}$.

\begin{lemma}\label{vectfield}
We have
$$Q_\e^{\circ q}(w) = w + 2i\pi q\e w(1-w^{q}) + \e R_\e(w),$$
with $R_\e\to 0$ uniformly on every compact subset of $\C$ as $\e\to 0$.
\end{lemma}

\proof
Since
$$P_{p/q+\e}^{\circ q} (z) = e^{2i\pi q\e}z+Az^{q+1}+{\cal O}(\e z^2),$$
we have
\begin{eqnarray*}
\frac{1}{z_\e}P_{p/q+\e}^{\circ q} (z_\e w) & = &
e^{2i\pi q\e}w + Az_\e^qw^{q+1} + {\cal O}(\e z_\e w^2)  \\
& = & w + 2i\pi q\e(w-w^{q+1}) + {\cal O}(\e^{1+1/q}w^2)
+ {\cal O}(\e^2 w).
\end{eqnarray*}
\qedprop

\begin{figure}[htbp]
\centerline{ \fbox{\includegraphics{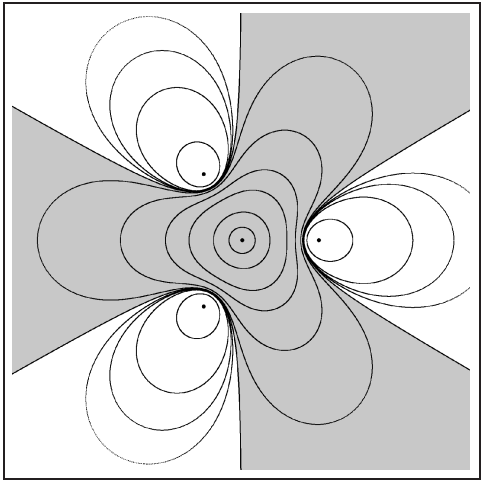}} } \caption{Some
trajectories of the real flow of the vector field $2i\pi
qw(1-w^q)\frac{\partial}{\partial w}$ for $q=3$.\label{field}}
\end{figure}

Figure \ref{field} shows some trajectories of the real flow of the
vector field $2i\pi q w(1-w^q)\frac{\partial}{\partial w}$ for
$q=3$. The origin is a center and its basin $\Omega$ is colored light grey.

Let us now define
$$Y_\e = \frac{1}{z_\e}Y\left(\frac{p}{q}+\e\right).$$
The set $Y_\e$ contains $1$ and we have
$$\log \rad(\C\setminus Y(p/q+\e)) = \log \rad(Y_\e) + \log |z_\e|.$$
Thus, we must show that
$$\limsup_{\e\to 0,~\e\in \R} \log \rad(\C\setminus Y_\e)\leq 0.$$

Set $\overline{Y_\e}=Y_\e\cup\{\infty\}$. This set is compact in $\P$.
Without loss of generality,
extracting a subsequence if necessary, we may assume that it converges
for the Hausdorff topology on compact subsets of $\P$ to some
limit $\overline{Y_0}$ as $\e\to 0$.
We define $Y_0=\overline Y_0\setminus \{\infty\}$.
Each $\overline{Y_\e}$ is connected and
contains $1$ and $\infty$. Passing to the limit, we see that $\overline{Y_0}$
is also connected and contains $1$ and $\infty$. Moreover, $Q_\e$
converges uniformly on compact subsets of $\C$ to the rotation $w\mapsto
e^{2i\pi p/q} w$. Since $Q_\e(Y_\e)=Y_\e$, we see that $Y_0$ is invariant
under this rotation.
Note that $Q_\e^{\circ q}(Y_\e)\subset Y_\e$ and
$$Q_\e^{\circ q}(w) = w + 2i\pi q\e w(1-w^q) + \e R_\e(w)$$
with $R_\e\to 0$ uniformly on compact subsets of $\C$ as $\e\to
0$. It follows that $Y_0$ is forward invariant under the real flow
of the vector field $2i\pi q w(1-w^q)\frac{\partial}{\partial w}$. Consider
the map $\phi:w\mapsto \zeta=w^q/(w^q-1)$. It is the composition of $w\mapsto
w^q$, (which identifies the quotient of $\P$ under the rotation of angle $1/q$
with $\P$), with a Möbius transformation fixing $0$, sending $1$ to $\infty$,
and $\infty$ to $1$. It sends the above vector field to the circular vector
field $(2\pi q^2) i\zeta \frac{\partial}{\partial \zeta}.$ It follows that
$Y_0$ contains the set $\phi^{-1}(\C\setminus\D)$. 
%We have already mentioned that the conformal radius of $\Omega$ is $1$.
Thus, we have
$$\limsup_{\e\to 0,~\e\in \R}\log\rad(\C\setminus Y_\e) \leq
\log \rad(\phi^{-1}(\D))= 0.$$ The proof of lemma \ref{sizeya} is
completed.

\section{Yoccoz's renormalization techniques.\label{rensec}}

In this section, we present the techniques of renormalization
developed by Yoccoz \cite{y}. We will follow the presentation
given by Pérez-Marco \cite{pm}. 

\remark There will be many constants in the discussion. Their
sharp value is not important for the application we will make here, so
we did not try to optimize them. Moreover, in many estimates where
$C \d$ appears, it can be weakened to $\epsilon(\d)$, where $\e(x)
\underset{x\to 0}{\tend} 0$, while still applying to our proof.

\subsection{Renormalization principle.\label{secyoc}}

Here, we recall what Pérez-Marco writes in \cite{pm} section~III,
adapting it to the setting of maps which are close to translations.

We denote by $T$ the translation $Z\mapsto Z+1$, by $S(\a)$ the
space of univalent mappings $F:\H\to \C$ such that $F\circ
T=T\circ F$ and such that $F(Z)-Z\to \a$ as $\im(Z)\to +\infty$.
This space is compact for the topology of uniform convergence on
compact subsets of $\H$.

Given $\d>0$, we denote by $S_\d(\a)$ the space of maps $F\in
S(\a)$ such that
\begin{equation}\label{fundest}
(\forall Z\in \H)\qquad |F(Z)-Z-\a|\leq \d\a \quad\text{and}\quad
|F'(Z)-1|\leq \d.
\end{equation}
Such a function $F$ extends continuously to $\H\cup \R$.

\vskip.2cm \noindent{\bf Step 1.}~ Assume $F\in S_\d(\a)$ and
define $l=i\R$ and $l'=[0,F(0)]$. If $\d$ is sufficiently small
(for example $\d<1/10$), $l\cup l'\cup F(l)$ bounds an open strip
${\cal U}$ in $\C$. Gluing the curves $l$ and $F(l)$ in the
boundary of $\overline {\cal U}$ via $F$, we obtain a surface ${\cal
V}$, whose remaining boundary corresponds to the segment
$l'$. Its interior is a Riemann surface for the complex 
structure inherited from $\overline{{\cal U}}$ (the gluing is
analytic). It is biholomorphic to the punctured disk $\D^*$.
Lifting via $Z\mapsto z=e^{2i\pi Z}$, we get an injective
holomorphic map $L:{\cal U}\to \H$ which extends continuously to
$\overline{\cal U}$ and such that
$$(\forall Z\in l)\qquad L(F(Z))=L(Z)+1.$$
We normalize $L$ by requiring $L(0)=0$.

\begin{proposition}\label{proxlem}
For all $\d\in ]0,1/10[$, all $\a\in ]0,1[$, all $F\in S_\d(\a)$, and
all $Z\in \overline {\cal U}$,
\begin{equation}\label{refyoc}
\im(Z)-2\d<\a\im(L(Z))<\im(Z)+2\d.
\end{equation}
\end{proposition}

\begin{proposition}\label{extlem}
Under the same assumptions, the map $L$ extends to a univalent map on
%${\cal W} \cup \overline{\cal U}$, with
$${\cal W}=\overline{\cal U} \cup\{Z\in \C~;~-1\leq \re(Z)\leq
  0~\text{and}~\im(Z)\geq 4\d\},$$
%% that conjugates $G$ to $T$ (more precisely, $L \circ G (z)= T \circ
%% L(z)$ for all $z$ such that both $z$ and $G(z)$ belong to ${\cal W}$).
%% Et même mieux !
\end{proposition}

From now on, $L$ will refer to this extension. The definition of
$\cal W$ is so that any point $Z\in W$ is
eventually mapped to $\cal U$ under iteration of $F$: $F^k(Z) = Z' \in
\cal{U}$ for some $k\in\N$. Then, one defines $L(Z)= L(Z')-k$. In
particular, $L$ conjugates $F$ to the translation $T$.

\vskip.2cm \noindent{\bf Step 2.}~ Given $\d\in ]0,1/10[$ and
$F\in S(\a)$, we can define inductively a sequence of univalent
maps $(F_n)_{n\geq 0}$ such that $F_n\in S(\a_n)$.
The construction depends on the choice at each step of some
real number $t_n>0$.
We start with $F_0=F-a_0$ (where $a_0=\lfloor \a\rfloor$) and we
assume that $F_n$ is constructed. We choose $t_n$ such that the
fundamental estimates (\ref{fundest}) hold for $\im(Z)\geq t_n$
(which is always possible).
It follows that $G_n:Z\mapsto F_n(Z+it_n)-it_n$ belongs to
$S_\d(\a_n).$ For $G_n$, we construct ${\cal U}_n$, ${\cal W}_n$ and
$L_n$ as above. Let $H_n$ be defined on $L_n\{Z\in\overline{\cal
U}~;~\Im(z)> 4\d\}$ by $H_n(z) = L_n \circ T^{-1}\circ L_n^{-1}$.
Note that, by proposition~\ref{proxlem}, if $\im(Z)>6\d/\a_n$,
there exists an integer $k$ such that $Z-k$ belongs to $D$, the domain of
definition of $H_n$.

% NOTE IMPORTANTE :
% 2 = C_0
% 4 = C_1
% 6 = C_1+C_0
% 8 = C_1+2.C_0 = B
% 32 = 4.B

Then, $D+\Z$ contains the
half plane ``$\Im(Z)>6\d/\a_n$''. Moreover,
the map $H_n$ commutes with the
translation $T$ on the set of points in $L_n(i[0,+\infty[)$ whose
imaginary part is $>6\d/\a_n$. This set being analytically
removable, this
implies $H_n$ extends univalently to the upper half-plane
$\{Z\in \C~|~\im(Z)> 6\d/\a_n\}$.
Moreover, as $\im(Z)\to +\infty$, $H_n(Z)-Z\to-1/\a_n=-a_{n+1}-\a_{n+1}.$

We set
$${\cal W}'_n={\cal W}_n+it_n$$ 
and we define $K_n:{\cal W}'_n\to \C$ by
$$K_n(Z)=s\circ{L_n(Z-it_n)}-i\frac{6\d}{\a_n}$$
where $s(x+iy)=-x+iy$,
and
$F_{n+1}\in S(\a_{n+1})$ defined on $\H$ by
$$F_{n+1}=K_n\circ T^{-1}\circ K_n^{-1}-a_{n+1}.$$
Note that on ${\cal W}'_n\cap F_n^{-1}({\cal W}'_n)$,
$K_n$ conjugates $F_n$ to $T^{-1}$.

\vskip.2cm
\noindent{\bf Step 3.}~
Next, to a point $Z\in \H$, we associate a sequence $(Z_n)_{n\geq
0}$ as follows. We define $Z_0=Z$. If $d_n=\im(Z_n)\geq 4\d+t_n$,
we choose $Z'_n$ such that $Z_n-Z'_n\in \Z$ and $-1\leq
\re(Z'_n)<0$, and we define
$$Z_{n+1}=K_n(Z'_n).$$
The sequence $(Z_n)_{n\geq 0}$ may be finite or infinite. The
estimates of proposition \ref{proxlem} imply that for $n\geq 0$
such that $Z_{n+1}$ is defined,
$$\im(Z_n)-t_n-8\d\leq \a_n\im(Z_{n+1})\leq \im(Z_n) 
-t_n-4\d.$$
For $n_0 \geq 0$:
\begin{equation}
\sum_{n=0}^{{n_0}-1}\beta_{n-1}(t_n+4\d)
\leq d_0 -\beta_{{n_0}-1}d_{n_0}
\leq \sum_{n=0}^{{n_0}-1}\beta_{n-1}(t_n+8\d)
\end{equation}
Which implies
\begin{equation}\label{heightdn}
\sum_{n=0}^{{n_0}-1}\beta_{n-1}t_n \leq d_0 -\beta_{{n_0}-1}d_{n_0} \leq
32\d+\sum_{n=0}^{{n_0}-1}\beta_{n-1}t_n
\end{equation}
%(1+\beta_0+\beta_1+\cdots+\beta_{n-2})B
Indeed, $1+\beta_0+\cdots+\beta_{n-2} \leq 4$
since $\beta_{-1}=1$, $\beta_0 = \a_0\leq 1$ and, $\beta_{n+2}\leq
\beta_n/2$.
%As a consequence, for all ${n_0}\geq 0$, we have
%%$$\beta_{n_0-1}d_{n_0}+\sum_{n=0}^{n_0-1}\beta_{n-1}t_n\leq d_0 \leq
%%\beta_{{n_0}-1}d_{n_0}+\sum_{n=0}^{{n_0}-1}\beta_{n-1}(t_n+12\d).$$
%% So, we have, for all ${n_0}\geq 0$,

\begin{proposition}\label{orbinf}
If $Z\in \H$ and if there exists $m\geq 0$ such that $F^{\circ m}(Z)\notin
\H$, then the sequence $(Z_n)_{n\geq 0}$ is finite.
\end{proposition}
\proof Let $H_n$ be the half plane
defined by ``$\Im Z>t_n$''.
If $Z_n$ is defined, let $1+k_n$ (with $k_n \geq 0$) be the rank of the first
iterate of $Z_n$ under $F_n : \H \to \C$ that leaves $H_n$. Note that if
$k_n=0$, then $Z_{n+1}$ is not defined.
Now, if $Z_{n+1}$ is defined and $k_{n+1}>0$, this means that $Z_{n}-k_{n+1}$
is eventually mapped back to ${\cal U}_n$ by iteration of $F_n$, without
leaving $H_n$. Therefore (since $|F_n(Z)-(Z+\a_n)|<
\a_n/10$ on $H_n$), $$k_{n+1} \leq \frac{11}{10}\a_n k_n.$$
Since $\a_n\a_{n+1} \leq 1/2$ this implies $k_{n+2} \leq \frac{121}{200}
k_n$ whenever defined, from which the proposition follows.
\qedprop

We can now reformulate Theorem III.1.1 in \cite{pm} as follows.

\begin{proposition}\label{size}
Assume we can choose the sequence $(t_n)_{n\geq 0}$ so that the
$n$-th renormalization $F_n$ satisfies the fundamental estimates
(\ref{fundest}) when $\im(Z)>t_n$ and so that
$$\Phi=\sum_{n=0}^{+\infty} \beta_{n-1}t_n <+\infty.$$
Then $F$ is linearizable and its Siegel disk contains the
following upper half-plane:
$$\left\{Z\in \C~|~\im(Z)>\Phi+32\d\right\}.$$
\end{proposition}

\proof It is enough to prove that all point $Z$ in the half plane has
infinite orbit. By proposition~\ref{orbinf}, this follows from the
sequence $Z_n$ being infinite. Indeed, assume $Z_n$ is
defined. According to the previous computations,
\begin{eqnarray*}
\beta_{n-1}d_n & \geq & d_0-\sum_{k=0}^{n-1}\beta_{k-1} t_k
-(1+\cdots+\beta_{n-2})8\delta \\
%& \geq & d_0-\Phi-(1+\cdots+\beta_{n-2} + \beta_{n-1})B\delta\\
& = &
(d_0-\Phi-32\d)+\beta_{n-1} t_n + \sum_{k=n+1}^{+\infty}\beta_{k-1}t_k +
\\
& & \big(4-(1+\cdots+\beta_{n-1})\big)8\delta+\beta_{n-1}8\d
\\ 
& \geq & \beta_{n-1}(t_n+8\delta)
\end{eqnarray*} 
Therefore,
$d_n \geq t_n+8\delta$. Since $8>4$, this implies $Z_{n+1}$ is 
defined.\qedprop

Also, there is a correspondence between periodic orbits for $F$ and for $F_n$.
Given a map $F : \H \to \C$ that commutes with $T$, we will say
that $Z\in\C$ is periodic with rotation number $p/q$ when
$F^q(Z)=Z+p$. In this case, $p$ and $q$ need not to be coprime.

\begin{proposition}\label{percor}
Let $n_0\geq 0$. If $F_{n_0}$ has a fixed point with rotation
number $0/1$ and imaginary part $h_{n_0}$,
then $F$ has a periodic orbit with rotation number $p_{n_0}/q_{n_0}$
contained in the strip
$$\{Z\in \C;~H\leq \im(Z)\leq H+32\d\}
\quad\text{with}\quad
H = \beta_{n_0-1}h_{n_0} + \sum_{n=0}^{n_0-1}\beta_{n-1}t_n.$$
Reciprocally, if $F$ has a periodic orbit with rotation number
$p_{n_0}/q_{n_0}$ whose imaginary part $h_0$ satisfies $\ds h_0 >
\sum_{n=0}^{n_0-1}\beta_{n-1}t_n +32\d$, then $F_{n_0}$ has a fixed
point of rotation number $0/1$, and height $h_{n_0}$ satisfying
$$h_0-32\delta\leq \beta_{n_0-1}h_{n_0} +
\sum_{n=0}^{n_0-1}\beta_{n-1}t_n\leq h_0.$$
\end{proposition}

\proof
Same as in \cite{pm} annex 2.e.
\qedprop

In the previous proposition, the reader should be aware that
$F_{n_0}(Z) = z+k$ with $k\in\Z^*$ is not considered as a fixed point with
rotation number $0/1$.

\subsection{Proof of proposition \ref{proxlem}.\label{proofproxlem}}

To obtain inequality (\ref{refyoc}) we will control of the
distortion of quasiconformal maps as follows.
Since $F(Z)-Z-\a$ is
periodic of period $1$, we have
$$|F(Z)-Z-\a|\leq \d\a e^{-2\pi \im(Z)}
\quad\text{and}\quad |F'(Z)-1|\leq \d e^{-2\pi \im(Z)}.$$ Let $B$
be the half-band $\{Z\in H~|~0<\re(Z)<1\}$. Let $H:\overline B \to
\overline{{\cal U}}$ be the map defined by
\begin{equation}\label{defh}
H(X+iY) = i\a Y + X \big[F(i\a Y)-i\a Y\big].
\end{equation}
An elementary computation shows that $\|\overline\partial
H/\partial H\|_\infty<1$ and if we set
$$K_H=\frac{1+|\overline \partial H/\partial H|}
{1-|\overline \partial H/\partial H|},$$
One computes that
\begin{eqnarray*}
\big|\partial H -\a\big| & \leq & \a \delta e^{-2\pi\a Y}\\
\big|\overline{\partial} H | & \leq & \a \delta e^{-2\pi\a Y}
\end{eqnarray*}
And therefore\footnote{A quick majoration yields a $4$, having a $2$
requires more care.}
\[K_H(X+iY)\leq \frac{1}{1-2\d e^{-2\pi \a Y}}.\]
Then, using $\d<1/10$, we
have the inequality
$$K_H(X+iY)\leq 1+\frac{5}{2}\d e^{-2\pi \a Y}.$$
In particular, $H$ is a $(1+\frac{5}{2}\d)$-quasiconformal homeomorphism.
Moreover, by definition
$$\im(H(Z))-\alpha\d\leq
\a\im(Z) \leq \im(H(Z))+\alpha\d,$$
and thus for all $Z\in \overline {\cal U}$, since $\alpha<1$:
$$\im(Z)-\d \leq \a\im(H^{-1}(Z)) \leq \im(Z)+\d.$$
Since $L$ is conformal, the map $G=L\circ H$ is quasiconformal
with the same dilatation as $H$. Moreover,
$G(iY+1) = G(iY)+1$ and so, since the imaginary axis is
quasiconformally removable, $G$ extends to a quasiconformal homeomorphism
$\H\to \H$. We will show that for all $Z\in \H$, we have
$$\a\im(Z)-\d \leq \a\im(G(Z)) \leq \a\im(Z)+\d.$$
It follows that
$$\im(Z)-2\d\leq \a\im(H^{-1}(Z))-\d \leq \a\im(L(Z)) \leq
\a\im(H^{-1}(Z))+\d \leq \im(Z)+2\d.$$

\begin{lemma}
Assume $\psi:(\D,0)\to (\D,0)$ is a $K$-quasiconformal
homeomorphism. Then, for all $z\in \D$,
$$4^{1-K}|z|^K\leq |\psi(z)|\leq 4^{1-1/K}|z|^{1/K}.$$
\end{lemma}

\proof To prove the upper bound, note that $\psi$ sends the
annulus $\D\setminus [0,z]$ to an annulus separating $0$ and
$\psi(z)$ from $S^1$. The modulus is divided by at most $K$. So,
$$|\psi(z)|\leq \mu^{-1}\left(\frac{\mu(|z|)}{K}\right),$$
where, for $r\in ]0,1[$, $\mu(r)$ is the modulus of the annulus
$\D\setminus [0,r]$ (it is a decreasing function). The estimate
$$\mu^{-1}\left(\frac{\mu(r)}{K}\right)\leq4^{1-1/K}r^{1/K}$$
can be found in \cite{avv} corollary 5.44.

The lower bound is obtained by applying the upper bound to
$\psi^{-1}$ which is $K$-quasiconformal. \qedprop

\begin{lemma}
If $\Psi:\H\to \H$ is a $K$-quasiconformal homeomorphism such that
$\Psi\circ T=T\circ \Psi$, then
$$\frac{1}{K}\im(Z)-\frac{K-1}{2\pi K}\log 4 \leq \im (\Psi(Z)) \leq
K\im(Z)+\frac{K-1}{2\pi}\log 4.$$
\end{lemma}

\proof $\Psi$ is the lift, via $Z\mapsto z=e^{2i\pi Z}$, of a
$K$-quasiconformal homeomorphism $\psi:(\D,0)\to (\D,0)$ as in the
previous lemma. \qedprop

We now come to the control of the quasiconformal homeomorphism
$G$.

\begin{lemma}
Let $\epsilon$ and $\eta$ be any two positive real numbers.
Assume $G:\H\to \H$ is a $(1+\e)$-quasiconformal homeomorphism
such that $G\circ T=T\circ G$ and
$$K_G(X+iY) \leq 1+\e e^{-\eta Y}.$$
Then,
$$\im(Z)-\frac{\e}{\eta} -\frac{\e}{2\pi (1+\e)}\log 4 \leq
\im(G(Z)) \leq\im(Z)+\frac{\e}{\eta} +\frac{\e}{2\pi}\log 4,$$
which yields
\[\big|\im(G(Z)) - \im(Z)\big| \leq \frac{\e}{\eta}
+\frac{\e}{2\pi}\log 4.\]
\end{lemma}

\proof We can write $G=G_2\circ G_1$ with
$$G_1(X+iY) = X+i\frac{1}{1+\e} \left(Y-\frac{\e}{\eta}e^{-\eta
Y}+\frac{\e}{\eta}\right).$$ An elementary computation shows that
$$K_{G_1}(X+iY) = \frac{1+\e}{1+\e e^{-\eta Y}}
\quad\text{and}\quad \im(G_1(Z))\leq
\frac{1}{1+\e}\left(\im(Z)+\frac{\e}{\eta}\right).$$ So, we can
apply the previous lemma to $G_2$ with $K=1+\e$, which yields the
upper bound for $\im(G(Z))$.

To get the lower bound, we use the same argument, writing
$G=G_4\circ G_3$ with
$$G_3(X+iY) = X+i(1+\e)\left(Y+\frac{1}{\eta}\log \frac{1+\e e^{-\eta
Y}}{1+\e}\right).$$ We have
$$K_{G_3}(X+iY) = \frac{1+\e}{1+\e e^{-\eta Y}}
\quad\text{and}\quad \im(G_3(Z))\geq
(1+\e)\left(\im(Z)-\frac{\e}{\eta}\right).$$
\qedprop

To conclude the proof of the proposition, we apply the previous lemma to
$\e=\frac{5}{2}\d$ and $\eta= 2\pi \a$. Using $\a<1$, we
have
$$\frac{\e}{\eta}+\frac{\e}{2\pi}\log 4 =
\frac{5\d}{4\pi\a}\big(1+\a\log 4\big) \leq \frac{\d}{\a}.$$

\subsection{Controlling the height of renormalization.}

In this section, we determine an upper bound for the height $t$ above which
the fundamental estimates (\ref{fundest}) are satisfied. The first result is
due to Yoccoz (it easily follows from the compactness of $S(0)$
, but the interested reader can find sharper bounds in \cite{y}, in the
lemma of section~3.2, page~26).

\begin{proposition}\label{yocprop}
For all $\d\in ]0,1/10[$, there exists a constant $C_\d$ such that
for all $F\in S(\a)$,
$$\im(Z)\geq C_\d\quad\Longrightarrow\quad |F'(Z)-1|\leq\d$$
and
$$ \im(Z)\geq \frac{1}{2\pi} \log
\frac{1}{\a}+C_\d\quad\Longrightarrow\quad 
|F(Z)-Z-\a|\leq\d\a.$$
\end{proposition}

(Of course, $C_\delta \tend +\infty$ when $\delta \tend 0$.)
\remark In particular, $F$ can not have fixed points above
$\frac{1}{2\pi} \log \frac{1}{\a}$ plus some universal constant.

The next result is a slight generalization of a result of
Pérez-Marco.

\begin{proposition}\label{pmprop}
For all $\d\in ]0,1/10[$, there exists a constant $C_\d$ such that
the following holds. Assume $\im(Z_0)\in \H$, $\a\in ]0,1[$ and
$F\in S(\a)$ has no fixed point except possibly $Z_0$ and its
translates by an integer. If
$$\im(Z)\geq \im(Z_0)+\frac{1}{2\pi}\left(\log \log \frac{e}{\a} -
\log (1+2\pi\im(Z_0))\right)+C_\d$$
then
$$|F(Z)-Z-\a|\leq\d\a.$$
\end{proposition}

One can rewrite
$$\log \log \frac{e}{\a} -
\log (1+2\pi\im(Z_0)) = \log \frac{1 + \log (\a^{-1})}{1 +
2\pi \im(Z_0)}.$$
Thus for $\im(Z_0) < \log (\a^{-1})/2\pi$, this number is positive.
From this, and the remark following proposition~\ref{yocprop}, it
follows that we can take the same constants $C_\delta$ in
propositions~\ref{yocprop} and~\ref{pmprop}.

\remark It follows that if $F$ has no fixed point, the
fundamental estimates (\ref{fundest}) are satisfied as soon as 
$$\im(Z)\geq \frac{1}{2\pi}\log\log\frac{e}{\a}+C_\d.$$
This result is due to Pérez-Marco \cite{pm}. This is the
form we will use in section~\ref{proofliminf1}.

\remark If $\ds \im(Z_0)\geq\frac{1}{2} \cdot
\frac{1}{2\pi}\log \frac{1}{\a}$, it follows from the two propositions
and an elementary computation that
the fundamental estimates (\ref{fundest}) are satisfied as soon as
$$\im(Z)\geq \im(Z_0)+1+C_\d.$$
This is the form\footnote{The assumption
$\im(Z_0)\geq\frac{1}{2} \cdot \frac{1}{2\pi}\log \frac{1}{\a}$ can be
replaced by $\im(Z_0)\geq \mu \cdot \frac{1}{2\pi}\log \frac{1}{\a}$
with $\mu \in ]0,1[$, giving the condition $\im(Z) \geq \im(Z_0) +
\log (\mu^{-1})/2\pi + C_\d$.} we will use in
section~\ref{cremerparabo}. 

\vskip.2cm

\noindent\textbf{Proof of proposition \ref{pmprop}.}
Without loss of generality, we may assume that
$$\im(Z_0)<\frac{1}{2\pi}\log \frac{1}{\a}$$
since otherwise, the result follows from proposition \ref{yocprop}. Let us set
$r=e^{-2\pi \im(Z_0)}$ if $F$ has a fixed point at $Z_0$ and $r=1$
if $F$ has no fixed point. Then, $\a<r$.

Let us now define $u(Z)=F(Z)-Z$. Since $u$ is $\Z$-periodic, there
exists a function $g:\D^*\to \C$ such that $u(Z)=g(e^{2i\pi Z})$.
The map $g$ extends holomorphically at $0$ by $g(0)=\a$. 
We need now to find an upper bound on $|z|$
which ensures that $|g(z)-\a|<\a \d$. By compactness of $S(0)$, we
can find a (universal) radius $r_0<1$ 
such that on $B(0,r_0)$, g takes its values in $B(0,e)$. Moreover,
if $F$ has a fixed point at $Z_0$, we define $\zeta_0=e^{2i\pi
Z_0}$. Then $g(\zeta_0)=0$ and $g$ does not vanish in $\D\setminus
\{\zeta_0\}$. If $F$ has no fixed point, $g$ does not vanish in
$\D$. In both cases, the map $g:B(0,r_0)\setminus \{\zeta_0\}\to
B(0,e)\setminus \{0\}$ is contracting for the hyperbolic metrics.

The coefficient of the hyperbolic metrics of $B(0,e)\setminus \{0\}$
at the point $\a$ is equal to $1/(\a \log(e/\a))$, so at first
approximation, points at hyperbolic distance of order
$\d/\log (e/\a)$ should be at Euclidean distance of order $\d\a$.
The lemma below makes a rigorous statement.
%% \begin{lemma}\label{leminutile}
%%   For all $\d \in ]0,1/10[$, for all $\a\in]0,1[$, for all $u \in \C$
%%   with $|u|=1$, 
%%   \[\frac{1}{2}\cdot\frac{\d}{\log(e/\a)} \leq d_{B(0,e)\setminus\{0\}}
%%   (\a,\a(1+\d u)) \leq 2\frac{\d}{\log(e/\a)}\]
%% \end{lemma}
%% \proof A uniformization of the universal cover of
%% $B(0,e)\setminus\{0\}$ is given by 
%% $z\in H \mapsto e\cdot \exp(
%% -\log(e/\a) \cdot\phi(z))$
%% where $H$ is the right half plane ``$\re z>0$''.
%% Therefore the hyperbolic distance from $\a$ to $\a (1+\d u)$ in
%% $B(0,e)\setminus\{0\}$ is equal to the hyperbolic distance in $H$ from
%% $1$ to $z = 1 + b$, where $b=\frac{\log(1+\d u)}{\log(e/\a)}$. Under
%% the further change of variable $v = \frac{1-z}{1+z} = \frac{b}{2+b}$,
%% this distance is equal to $2\on{argth}|v|$. Let
%% $x=\d/\log(e/a)$. We then use convexity
%% inequalities to bound $|b|/x$, then $|v|/x$ and finally
%% $2\on{argth}|v|/x$.
%% %% $0.9 x \leq |b| \leq 1.1 x$
%% %% $1.8 \leq |2+b|\leq 2.2$
%% %% $0.4 x \leq |v| \leq 0.62 x$
%% %% $|v| \leq 0.056$
%% %% $0.8 x \leq 2\on{argth} |v| \leq 1.25 x$
%% \qedlem
%% The lower inequality yields
%% \begin{corollary}\label{corda}
%% \textcolor{black}{
%% $(\forall \d\in ]0,1/10[)$, $(\forall \a\in ]0,1[)$,
%% $$d_{B(0,e)\setminus \{0\}}(\a,z)\leq\frac{\d}{2\log e/\a}
%% \quad\Longrightarrow\quad |z-\a|\leq \d\a.$$}
%% \end{corollary}
\begin{lemma}\label{lemda}
$(\forall \d\in ]0,1/10[)$, $(\forall \a\in ]0,1[)$,
$$d_{B(0,e)\setminus \{0\}}(\a,z)\leq\frac{\d}{2\log e/\a}
\quad\Longrightarrow\quad |z-\a|\leq \d\a.$$
\end{lemma}
\proof For $x<\a$, let $\rho(x)$ be the infimum of the
coefficient of the hyperbolic metric on the Euclidean circle of center
$\a$ and radius $x$.
If $|z-\a|>\d\a$, then
the hyperbolic geodesic in $B(0,e)\setminus\{0\}$ from $\a$ to
$z$ is longer than \[\ds \int_0^{\d \a}\rho(x)dx.\]
Let us introduce the function
\[\renewcommand{\arraystretch}{1.5}\left\{\begin{array}{lll}
f(x) = \frac{1}{x\log e/x} & & 0<x \leq 1\\
f(x) = 1 & & 1\leq x<e
\end{array}\right.\]
Then $f$ is decreasing, and $\rho(x)=f(x+\a)$. Moreover,
$f$ is $C^1$ and convex, and therefore above its
tangents. Therefore
\begin{eqnarray*}
d_{B(0,e)\setminus \{0\}}(\a,z) & \geq & \int_{\a}^{\a+\d\a} f(x)dx \\
& \geq & \int_{\a}^{\a+\d\a} \big( f(\a) + (x-\a) f'(\a)\big) dx\\
& = & \frac{\d}{\log e/\a} \left(1- \frac{\d}{2}\Big(1- \frac{1}{\log
 (e/\a)}\Big)\right) \geq c \frac{\d}{\log e/a}
\end{eqnarray*}
with $c = 19/20 >1/2$.
%% \begin{eqnarray*}
%%  d_{B(0,e)\setminus
%%  \{0\}}(\a,z)
%%  & > \ds\int_{\a}^{\a(1+\d)} \frac{1}{t\log e/t}dt
%%  & = -\log \left(1-\frac{\log (1+\d)}{\log e/\a}\right) \\
%%  && \geq  \frac{\log (1+\d)}{\log e/\a}\geq \frac{\d}{2\log e/\a}.
%% \end{eqnarray*}
\qedlem

The next lemma is also motivated by a hyperbolic metrics coefficient
computation.
\begin{lemma}
$(\forall r_0<1)$, $(\exists
\gamma>0)$, $(\forall \d\in ]0,1/10[)$, if $0<\a<r\leq 1$, then
$$|z|\leq \gamma \d r \frac{\log e/r}{\log e/\a}
\quad\Longrightarrow\quad d_{B(0,r_0)\setminus\{r\}}(0,z)\leq
\frac{\d}{2\log e/\a}.$$
\end{lemma}

\proof
First case: $r\geq r_0/2$.\\
When $|z|\leq \d r_0$, then
$$d_{B(0,r_0)\setminus\{r\}}(0,z)\leq
d_{B(0,r_0/2)}(0,z)= \log \frac{1+2|z|/r_0}{1-2|z|/r_0} \leq
\frac{5 |z|}{r_0}.$$ Thus, when $r\geq r_0/2$, we can take any
$\gamma$ such that
$$\gamma\leq \min_{r\in [r_0/2,1]}\frac{r_0}{10 r\log e/r}
=\frac{1}{10\log e/r_0}.$$
Second case: $r< r_0/2$.\\
We first solve the problem when $r_0=1$.
Let $\rho(z)|dz|$ be the element of hyperbolic metric on
$\D\setminus\{r\}$. A computation gives
\[\rho(z) = \frac{1-r^2}{|1-rz| \cdot |z-r| \cdot
  \log\Big(\frac{|1-rz|}{|z-r|}\Big)}.\] 
A majoration gives, for $|z|<r/10$, $\rho(z)<10/(9r\log |s|^{-1})$ with
$s=(z-r)/(1-rz)$. Then, $|s|< 11r/(10+r^2)< 11r/10$.
Thus
\[\forall r\in]0,1/2[,\ \forall z\ \text{with}\ |z|\leq\frac{r}{10},\
\rho(z) \leq \frac{12}{r\log e/r}.\]
Therefore, for $r_0=1$, we can take $\gamma=\gamma_1$, with
\[\gamma_1 = 12.\]
For $r_0\in]0,1[$, we rescale the problem by the factor $1/r_0$,
and according to what we did above, a sufficient condition on
$z$ is that
\[\left|\frac{z}{r_0}\right| < \gamma_1 \d \frac{r}{r_0} \frac{\log
e r_0/r}{\log e/\a},\]
Then, using $r<r_0/2$, we can take
\[\gamma \leq \gamma_1 \frac{\log 2e}{\log 2e  + \log r_0^{-1}}.\]
%% and define
%% $$\rho = \gamma \d r \frac{\log e/r}{\log e/\a}.$$
%% Note that when $\gamma \d \leq 1/2$, then $\rho\leq r/2$. If $|z|\leq
%% \rho$, then,
%% $$d_{B(0,r_0)\setminus\{r\}}(0,z)\leq d_{B(0,r_0)\setminus\{r\}}(0,\rho) =
%% \log \frac{\ds \log\frac{1-r\rho/r_0}{r/r_0-\rho}}{\ds
%% \log\frac{r_0}{r}}.$$
%% Since $r_0<1$, $r/r_0\leq 1/2$, $\rho/r\leq
%% 1/2$ and $\log u\leq u-1$, an elementary computation gives
%% $$  d_{B(0,r_0)\setminus\{r\}}(0,z)\leq \frac{2\rho}{r\log r_0/r}.$$
%% Thus, when $r\leq r_0/2$, we can take any $\gamma\leq 5$ such
%% that
%% $$\gamma\leq \min_{r\in [0,r_0/2]} \frac{\log r_0/r}{4\log e/r}.$$
\qedlem

The two previous lemmas show that there
exists $\gamma>0$ such that for all $\d\in ]0,1/10[$, 
$$|z|\leq \gamma \d r\frac{\log e/r}{\log e/\a}
\quad\Longrightarrow\quad |g(z)-\a|\leq \d\a.$$
As a consequence,
$$\im(Z)\geq \frac{1}{2\pi}\left(\log \frac{1}{\gamma\d}+\log
  \frac{1}{r}+\log \frac{\log e/\a}{\log e/r}\right)
\quad\Longrightarrow\quad |F(Z)-Z-\a|\leq \d\a.$$
\qedprop

\section{Proof of inequality (\ref{liminf}) in most cases.\label{proofliminf1}}

We will use the following fact several times.
Assume $\a'\in ]0,1[$ tends to $\a\in ]0,1[$ and $F_{\a'} \in
S(\a')$ tends to the translation $T_\a:Z\mapsto Z+\a$ uniformly on
every compact subset of $\H$. Then, the convergence is uniform on
every upper half-plane of the form 
``$\im(Z) \geq t > 0$'', and $F'_{\a'} \tend 1$ uniformly on
these half-planes. Therefore,
given $\d\in ]0,1/10[$ and $t>0$, if $F_{\a'}$ is sufficienlty close to
$T_\a$,  the map $G_{\a'}:Z\mapsto F_{\a'}(Z+it)-it$ belongs to
$S_\d(\a')$ (it is important that $\a\neq 0$).
For $G_{\a'}$ we can construct ${\cal U}_{\a'}$, ${\cal W}_{\a'}$ and
$L_{\a'}$ as in section \ref{secyoc}. We then define
${\cal W}'_{\a'}={\cal W}_{\a'} + it$,
$K_{\a'}:{\cal W}'_{\a'}\to \C$ by
$$K_{\a'}(Z)=s\circ L_{\a'}(Z-it)-i\frac{6\d}{\a'}$$
and $F_{\a',1}\in S(\a'_1)$ by
$$F_{\a',1}= K_{\a'}\circ  T^{-1}\circ K_{\a'}^{-1} -\left\lfloor
\frac{1}{\a'}\right\rfloor,$$
where $s(Z) = -\overline{Z}$.

As $F_{\a'}$ tends to $T_\a$, $L_{\a'}$ tends to $Z\mapsto Z/\a$
uniformly on every compact subset of ${\cal W}$. Indeed, as in
section \ref{proofproxlem}, we can write $L=G\circ H^{-1}$ where
$H$ is defined by equation (\ref{defh}). Then, $H$ converges to
$Z\mapsto \a Z$ uniformly on $\overline B$ and $G:\overline \H\to
\overline \H$ is a $K$-quasiconformal homeomorphism such that
$G(0)=0$ and $G\circ T=T\circ G$. Moreover, $K\tend 1$ as
$F_{\a'}\tend T_\a$. Thus $G$ converges to the
identity uniformly on every compact subset of $\overline \H$.

%% As $F_{\a'}$ tends to $T_\a$, $L_{\a'}$ tends to $Z\mapsto Z/\a$
%% uniformly on ${\cal W}_{\a}
%% = \big\{Z\in\H\big|-1\leq \Re(Z) \leq \a\big\}$.
%% This follows from proposition~\ref{proxlem} applied to
%% $\delta \tend 0$.

It follows that $K_{\a'}$ tends to $Z\mapsto (s(Z)-it-i6\d)/\a$
uniformly on every compact subset of ${\cal W}'_{\a'}$ and
$F_{\a',1}$ tends to the translation
$Z\mapsto Z+\a_1$ uniformly on every compact subset of $\H$.

\subsection{Brjuno numbers.\label{secsiegel}}

Assume $\a\in ]0,1[$ is a Brjuno number and let $\phi_\a:\D\to
\Delta_\a$ be a linearizing parameterization. Note that
$|\phi_\a'(0)|=r(\a)$. For $\a'$ close to $\a$, let us define
$$f_{\a'}=\phi_\a^{-1}\circ P_{\a'}\circ \phi_\a$$
on $\phi_\a^{-1}(\Delta_\a\cap P_{\a'}^{-1}(\Delta_\a))$. Since
$P_\a(\Delta_\a)=\Delta_\a$ and $P_{\a'}\tend P_\a$ as $\a'\tend \a$,
we see that when $\a'\tend \a$, $f_\a$ converges uniformly on every
compact subset of $\D$ to the rotation of angle $\a$. Note that
when $\a'$ is a Brjuno number, $f_{\a'}$ has a Siegel disk of
radius $\rho(\a') \leq  r(\a')/r(\a).$ Indeed, the image of this
Siegel disk by $\phi_\a$ is contained in the Siegel disk of
$P_{\a'}$. Finally, let $F_{\a'}$ be the lift of
$f_{\a'}$ via $Z\mapsto e^{2i\pi Z}$ which satisfies
$|F(Z)-Z-\a'| \tend 0$ when $\Im(Z) \tend +\infty$.

Let us now fix $\eta>0$, $\d\in ]0,1/10[$ and $n_0\geq 1$. For
$n\geq 0$, we will define a sequence of heights $t'_n$ and a
sequence of maps $F_{\a',{n+1}}\in S(\a'_{n+1})$ as in section
\ref{secyoc}.

%% \begin{lemma}\label{choicet'n}
%% If $\a'\in \R\setminus \Q$ is sufficiently close to $\a$, we can
%% take
%% $$t'_0=\ldots=t'_{n_0} = \eta/(n_0+1).$$
%% \end{lemma}

%% \proof
%% As $\a'$ tends to $\a$, the domain of $F_{\a'}$ eventually contains the
%% half-plane $\{Z\in \C~|~\im(Z)>\eta/(n_0+1)\}$ and $F_{\a'}$ tends to the
%% translation
%% $Z\mapsto Z+\a$ uniformly this half-plane. As a consequence, if $\a'$ is
%% sufficiently close to $\a$, we can take $t'_0=\eta/(n_0+1)$.
%% Now, as $\a'$ tends to $\a$, $F_{\a',1}$
%% converges uniformly to the translation $Z\mapsto Z+\a_1$ on every compact
%% subset of $\H$. So, we can take $t_1=\eta/(n_0+1)$, and so on\ldots
%% \qedprop

According to the fact mentioned at the beginning of
section~\ref{proofliminf1}, and using induction on $n_0$, we know
that provided $\a'\in \R\setminus \Q$ is sufficiently
close to $\a$, we can take $$t'_0=\ldots=t'_{n_0}
= \eta/(n_0+1).$$

By proposition \ref{yocprop}, for $n\geq n_0+1$, we can take
$$t'_n=\frac{1}{2\pi}\log \frac{1}{\a'_n}+C_\d$$
for some constant $C_\d$ which only depends on $\d$.

It follows from proposition \ref{size} that if $\a'\in {\cal B}$ is
sufficiently close to $\a$, we have
\begin{eqnarray*}
\log \frac{r(\a)}{r(\a')}\ \ \leq\ \ \log \frac{1}{\rho(\a')}
& \leq & 2\pi \left(\sum_{n=0}^\infty \beta'_{n-1}
t'_n +32\d \right)\\
& \leq & \Phi(\a')-\Phi_{n_0}(\a')+2\pi(\eta +4\beta'_{n_0}C_\d +
32\d)
\end{eqnarray*}
(we used $\beta'_{n_0} + \beta'_{n_0+1} + \ldots \leq 4 \beta'_{n_0}$ which
follows from $\beta'_{k+1}\leq \beta'_k$ and $\beta'_{k+2} \leq \beta'_k
/2$).
Let us rewrite it
$$\Phi(\a') +\log r(\a') \geq \Phi_{n_0}(\a') + \log r(\a) -2\pi(\eta
+4\beta'_{n_0}C_\d + 32\d).$$
Letting $\a' \tend \a$ and using
$\Phi_{n_0}(\a') \tend \Phi_{n_0}(\a)$ and $\beta'_{n_0} \tend \beta_{n_0}$, 
$$\liminf_{\a'\to \a,~\a'\in {\cal B}} \Phi(\a') +\log r(\a') \geq
\Phi_{n_0}(\a) + \log r(\a) -2\pi(\eta +4\beta_{n_0}C_\d + 32\d).$$
Now, as $n_0 \tend +\infty$, $\Phi_{n_0}(\a) \tend \Phi(\a)$ and
$\beta_{n_0} \tend 0$. Thus
$$\liminf_{\a'\to \a,~\a'\in {\cal B}} \Phi(\a') +\log r(\a') \geq
\Phi(\a) + \log r(\a) -2\pi(\eta + 32\d).$$
Since this is valid for all $\eta>0$ and $\delta \in ]0,1/10[$, it
implies
$$\liminf_{\a'\to \a,~\a'\in {\cal B}} \Phi(\a') +\log r(\a') \geq
\Phi(\a) + \log r(\a) = \Upsilon(\a).$$
%$$\liminf_{\a'\to \a,~\a'\in {\cal B}} \Phi(\a')+\log r(\a') \geq
%\Phi(\a)+\log r(\a)-\e=\Upsilon(\a)-\e.$$

\subsection{Rational numbers.}

We consider a rational number $\a=p/q$ and a Brjuno number $\a'$
close to $p/q$.
Let us note $\a'_n$ and $\beta'_n$ the sequences associated to
$\a'$. According to the sign of $\e=\a'-p/q$, we
associated in section \ref{subsec_urat}
to $\a=p/q$ an integer $n_0 \in\N$, and finite sequences $\a_0$,
$\a_1$, \ldots, $\a_{n_0}=0$,
and $p_0/q_0$, $p_1/q_1$, \ldots, $p_{n_0}/q_{n_0} = p/q$
such that for all $k\leq n_0$, $\a'_k \tend \a_k$, $p'_k \tend p_k$ and
$q'_k \tend q_k$ when $\a' \tend \a$ on one side.

We will use the notations of section \ref{parexparabo}.
Let $z_\e$ be a point of the cycle ${\cal C}_{p/q}(\a')$.
To study the dynamics of $P_{p/q+\e}$ at the scale
of $z_\e$, we defined
$$Q_\e:w\mapsto \frac{1}{z_\e} P_{p/q+\e}(z_\e w).$$
Lemma \ref{vectfield} asserts that
\begin{equation}\label{eq_vfc}
Q_\e^{\circ q}(w) = w + 2i\pi q\e w(1-w^{q}) + \e R_\e(w),
\end{equation}
with $R_\e\tend 0$ uniformly on every compact subset of $\C$ as $\e\tend 0$.

Set $\phi(w) = \omega^q/(1-\omega^q)$
and $\Omega=\phi^{-1}(\D)$. It is the preimage by
$w \mapsto w^q$ of the half plane ``$\Re(z) <1/2$'' and is illustrated
as a gray set for $q=3$ in figure~\ref{field} page~\pageref{field}.
Let $\psi:\Omega\to \D$ be a holomorphic map satisfying
$\psi(w)^q=\phi(w)$. Then, $\psi(0)=0$, 
$|\psi'(0)|=1$ and $\psi$ is a conformal representation between
$\Omega$ and $\D$. It sends the vector field $2i\pi q
w(1-w^q)\frac{\partial}{\partial w}$ to the vector field $2i\pi
q\zeta \frac{\partial}{\partial \zeta}.$ We
define
$$f_\e =\psi\circ Q_\e\circ \psi^{-1}$$
on $\psi(\Omega \cap Q_\e^{-1}(\Omega))$. As $\e\tend 0$, $f_\e$
converges uniformly on every compact subset of $\D$ to the
rotation of angle $p/q$. Moreover
by (\ref{eq_vfc}) we see that when $\e\tend 0$,
$$f_\e^{\circ q}  (z) = z+2i\pi q\e z+ \e g_\e(z),$$
with $g_\e\tend 0$ uniformly  on every compact subset of
$\overline \D$.
Note that when $\a'=p/q+\e$ is a Brjuno
number, $f_\e$ has a Siegel disk of conformal radius
$$\rho(\e) \leq r(\a')/|z_\e|.$$
Let $F_\e$ be
the lift of $f_\e$ via $Z\mapsto e^{2i\pi Z}$ which satisfies
$|F_\e(Z)-Z-\a'| \tend 0$ when $\Im(Z) \tend +\infty$.
When $\e\tend 0$,
$$F_\e^{\circ q}\circ T^{-p}(Z) = Z + q\e + \e G_{\e}(Z)$$
with $G_\e\tend 0$ uniformly  on every compact subset of $\H$.

Let us fix $\d\in ]0,1/10[$ and $\eta>0$. For $n\geq 0$, we will
define a sequence of heights $t'_n$ and a sequence of maps
$F_{\e,{n+1}}\in S(\a'_{n+1})$.

As $\e$ tends to $0$, $F_\e$ converges uniformly to the
translation by $p/q$ on the upper half-plane $\{Z\in
\C~|~\im(Z)\geq \eta/(n_0+1)\}$. Moreover, for $n\leq n_0-1$, as $\e\tend
0$, $\a'_n\tend \a_n\neq 0$. Thus, if $\e$ is
sufficiently close to $0$, we can take
$$t'_0=t'_1=\ldots=t'_{n_0-1} = t \egaldef \eta/(n_0+1).$$
We will call ${\cal W}'_{\e,n}$ and
$K_{\e,n} : {\cal W}'_{\e,n} \to \C$
the objects corresponding to $W'_n$ and $K_n$ defined in
section~\ref{secyoc}. When $\e \tend 0$, the interior of ${\cal W}'_{\e,n}$
tends to the interior of a set
${\cal W}'_{0,n}$ which is the union of two half strips $``-1\leq \Re( Z)
\leq 0 \text{ and }
\Im (Z) \geq 4\d + t "$ and  $``
0\leq \Re (Z) \leq \a_n \text{ and }
\Im (Z) \geq t"$.
%% $${\cal W}'_{0,n} = \big\{Z\in \C\ \big|\ -1\leq \Re Z \leq 0 \text{ and }
%% \Im Z \geq 4\d + \eta/(n_0+1) \big\} \cup \big\{Z\in\C\ \big|\ 
%% 0\leq \Re Z \leq \a_n \text{ and }
%% \Im Z \geq \eta/(n_0+1)\big\}.$$}
For $n\leq n_0-1$,
as $\e$ tends to $0$, $K_{\e,n}$ tends to $Z\mapsto
(s(Z)-it-i6\d)/\a_n$ uniformly on every compact subset of
${\cal W}'_{0,n}$, where $s(Z) = -\overline{Z}$.

Now, when $\e\tend 0$, $F_{\e,n_0}$ converges uniformly to the
translation $Z\mapsto Z+\a_{n_0} = Z+0$, i.e., to the identity.

\begin{lemma}
If $\e$ is small enough, we can take
$t'_{n_0}=\eta/(n_0+1)$.
\end{lemma}

\proof
Let us now consider the map
$$\Psi_\e=K_{\e,n_0-1}\circ \ldots \circ K_{\e,0}.$$
Its set of definition eventually contains every compact subset
of the interior of $${\cal W}''=\big\{Z\in \C~;~
-\beta_{n_0-1}\leq (-1)^{n_0} \re(Z)\leq \beta_{n_0-2} ~{\rm
and}~\im(Z)\geq t'-2\d\beta_{n_0-2}\big\},$$
 with
$t'=(t+6\d)(1+\beta_1+\ldots+\beta_{n_0-2})$. On every of these compact
subsets, $\Psi_\e$ eventually conjugates
$F_\e^{\circ q}\circ T^{-p}$ to $F_{\e,n_0}$.

As $\e$ tends to $0$, $\Psi_\e$ converges to $Z\mapsto
(s^{n_0}(Z)-it')/\beta_{n_0-1}$, uniformly on every compact subset
of the interior of ${\cal W}''$.
Thus, since $s^{n_0}\circ\Psi_\e$ is holomorphic, the
derivative of $s^{n_0}\circ\Psi_\e$ converges to
$1/\beta_{n_0-1}$, uniformly on every compact subset of the
interior of ${\cal W}''$. Therefore $$F_{\e,n_0}(Z) = Z+ \frac{q
|\e|}{\beta_{n_0-1}} + \e H_\e(Z)$$ with $H_\e\tend 0$ uniformly  on
every compact subset of $\H$. Since
$\a'_{n_0}=q|\e|/\beta'_{n_0-1}
=q|\e|/\beta_{n_0-1}+{\cal O}(\e^2)$,
$|F_{\e,n_0}(Z) - Z -\a'_{n_0}| = \a'_{n_0} I_\e(Z)$
with $I_\e(Z) \tend 0$ uniformly on every compact subset of $\Psi_0 ({\cal W
}'')$. 
This set contains $``-1< \Re(Z) <1 \text{ and } \Im(Z)>0"$. 
Since $F_{\e,n_0}$ commutes with $T$, this implies that 
$|F_{\e,n_0}(Z) - Z -\a'_{n_0}| = \a'_{n_0} I_\e(Z)$
with $I_\e(Z) \tend 0$ uniformly on every compact subset of $\H$.
As a consequence
$|F'_{\e,n_0}(Z) - 1| \tend 0$ uniformly on every compact subset of
$\H$.
 \qedprop

Finally, for $n\geq n_0+1$, we can take
$$t'_n=\frac{1}{2\pi}\log \frac{1}{\a'_n}+C_\d$$
where $C_\d$ is the constant in proposition \ref{yocprop}. So, if
$\e$ is sufficiently small, we have
$$\log \frac{|z_\e|}{r(\a')} \leq 2\pi\left(\sum_{n=0}^\infty
\beta'_{n-1} t'_n+32\d\right)\leq
\Phi(\a')-\Phi_{n_0}(\a')+2\pi(\eta+4\beta'_{n_0}C_\d+32\d).$$
Reordering the terms, we obtain
$$\Phi(\a')+\log r(\a') \geq \log |z_\e|+\Phi_{n_0}(\a')
-2\pi(\eta+4\beta'_{n_0}C_\d+32\d).$$ As $\e\tend 0$, $\log
|z_\e|+\Phi_{n_0}(\a')$ tends to $\Upsilon(p/q)$ and
$\beta'_{n_0}$ tends to $0$. We therefore have (see lemma~\ref{moduleze})
$$\liminf_{\a'\to p/q,~\a'\in {\cal B}}\Phi(\a')+\log r(\a')
\geq \Upsilon\left(\frac{p}{q}\right)-2\pi (\eta+32\d)$$ and the
proof of inequality (\ref{liminf}) at rational numbers is
completed since $\eta$ and $\d$ can be chosen arbitrarily small.

\subsection{\Cremer numbers whose Pérez-Marco sum converges.\label{loglogsec}}

It is possible to give a proof that works for all \Cremer numbers
at the same time, but for clarity, we prefer to study two cases
(which overlap) separately. Here, we will assume
$\a$ is a \Cremer number such that
$$\sum_{n=0}^\infty \beta_{n-1}\log \log \frac{e}{\a_n}<\infty.$$
We will call this sum the Pérez-Marco sum, since it
was introduced by Pérez-Marco in \cite{pm}.
There, he proves that, under this condition, every
germ that fixes $0$ with derivative $e^{2i\pi \a}$ is linearizable
or has small cycles.

%$$\sup_n \beta_{n-1}\log \frac{1}{\a_n} = \infty.$$
%$$\sum_{n=0}^\infty \beta_{n-1}\log \log \frac{e}{\a_n}<\infty.$$

Let us fix $\eta>0$, $\d\in ]0,1/10[$ and $n_0\geq 1$. For
$n_1\geq n_0$, we set
$$d_{n_1}(\a')=d(0,X_{n_1}(\a'))$$ (see definition \ref{def_X}
for $X_n$). Since a Cremer point of a polynomial is 
accumulated by periodic points, and because we defined $X_{n_1}(\a)$
as the set of all periodic points of period $\leq q_{n_1}$ except $0$,
we have $d_{n_1}(\a) \tend 0$ when $n_1 \tend
+\infty$. Thus, provided $n_1$ is big enough, we see that for all $\a'$ close
enough to $\a$, $F_{\a'}$ is injective on $B(0,d_{n_1}(\a'))$.
Let $F_{\a'}\in S(\a')$ be the lift of $P_{\a'}$ via $Z\mapsto
d_{n_1}(\a')e^{2i\pi Z}$. This amounts to restrict
the polynomial $P_{\a'}$ to the disk $B(0,d_{n_1}(\a'))$ where
there are no periodic cycle of period less than or equal to
$q_{n_1}$, except $0$. Note that when $\a'$ is a Brjuno number,
this restriction has a Siegel disk of conformal
radius $\leq r(\a')$.

For $n\geq 0$, we will define a sequence of heights $t'_n$ and a
sequence of maps $F_{\a',{n+1}}\in S(\a'_{n+1})$.

\begin{lemma}
If $n_1$ is sufficiently large and $\a'$ is sufficiently close to $\a$, we can
take $$t'_0=t'_1=\ldots=t'_{n_0}=\eta/(n_0+1).$$
\end{lemma}

\proof Let us choose $\e$ sufficiently small so that $\a'_0\neq
0$, $\ldots$, $\a'_{n_0}\neq 0$ for all $\a'\in [\a-\e,\a+\e]$. As
$n_1\tend \infty$, $(\a',Z)\mapsto F_{\a'}(Z)-Z-\a'$ converges
uniformly to $0$ on $[\a-\e,\a+\e]\times \{Z\in \C~|~\im(Z)\geq
\eta/(n_0+1)\}$. If $n_1$ is sufficiently large, we can therefore
take $t'_0=t'_1=\ldots=t'_{n_0}=\eta/(n_0+1)$. \qedprop

By construction, the maps $F_{\a'}$ have no periodic cycle of period less than
or equal to $q_{n_1}$. So, by proposition~\ref{percor}, for
$n\leq n_1$, the renormalizations $F_{\a',n}$ have no fixed point
in $\H$. Thus, by proposition \ref{pmprop}, we can take
$$t'_{n_0+1}=\frac{1}{2\pi}\log\log\frac{e}{\a_{n_0+1}}+C_\d
\quad\ldots\quad
t'_{n_1}=\frac{1}{2\pi}\log\log\frac{e}{\a_{n_1}}+C_\d$$
for some constant $C_\d$ which only depends on $\d$.
Finally, by proposition \ref{yocprop}, for $n\geq n_1+1$, we can take
$$t'_n=\frac{1}{2\pi}\log\frac{1}{\a'_n}+C_\d.$$

Now, proposition \ref{size} yields
$$\frac{1}{2\pi}\log \frac{d_{n_1}(\a')}{r(\a')} \leq\sum_{n=0}^{\infty}
  \beta'_{n-1} t'_n +32\d.$$
Using the value of $t'_n$ chosen above, we get
%\begin{equation}\label{limitcremer1}
%\begin{array}{rcl}
\begin{eqnarray*}
\Phi(\a')+\log r(\a') &\geq &\Phi_{n_1}(\a')+\log d_{n_1}(\a')-
\ds \sum_{n=n_0+1}^{n_1}\beta'_{n-1}\log\log\frac{e}{\a'_n} \\
&& -2\pi(\eta+4\beta'_{n_0}C_\d+32\d).
\end{eqnarray*}
%\end{array}
%\end{equation}
Let $\a'$ tend to $\a$:
\begin{eqnarray*}
\liminf_{\a' \to \a,~\a'\in {\cal B}}\Phi(\a')+\log r(\a') & \geq
& \Phi_{n_1}(\a)+\log d_{n_1}(\a)- 
\ds \sum_{n=n_0+1}^{n_1}\beta_{n-1}\log\log\frac{e}{\a_n} \\
 & & -2\pi(\eta+4\beta_{n_0}C_\d+32\d).
\end{eqnarray*}
Let $n_1$ tend to $+\infty$. Remind that $d_{n_1}(\a) \sim
r_{n_1}(\a)$, and by definition $\ds \Upsilon(\a) = \liminf_{n_1 \tend
+\infty} \Phi_{n_1}(\a) + r_{n_1}(\a)$. Thus,
\begin{eqnarray*}
\liminf_{\a' \to \a,~\a'\in {\cal B}}\Phi(\a')+\log r(\a') & \geq
& \Upsilon(\a)- 
\ds \sum_{n=n_0+1}^{+\infty}\beta_{n-1}\log\log\frac{e}{\a_n}
\\ & &  -2\pi(\eta+4\beta_{n_0}C_\d+32\d).
\end{eqnarray*}
Let $n_0$ tend to $+\infty$. Since $\beta_{n_0} \tend 0$ and
the Pérez-Marco sum of $\a$ was assumed to be convergent, we have
\begin{eqnarray*}
\liminf_{\a' \to \a,~\a'\in {\cal B}}\Phi(\a')+\log r(\a') & \geq
& \Upsilon(\a) -2\pi(\eta+32\d).
\end{eqnarray*}
Since this is valid for arbitrarily small $\eta$ and $\delta$, this
concludes the proof for the case when the Pérez-Marco sum of $\a$ converges.

%% Given $\e>0$, if $\eta$  small enough, we have $2\pi\eta<\e/6$, if
%% $\d$ is small enough, we have $64\pi \d<\e/6$ and if $n_0$ is
%% large enough, $8\pi \beta_{n_0}C_\d<\e/6$ and
%% $$\sum_{n=n_0+1}^{\infty}\beta_{n-1}\log\log\frac{e}{\a_n}<\e/6$$
%% (this is possible since the series converges).
%% Now, close to $0$, $P_\a$ is close to the rotation of angle $\a$, and since
%% $d_n(\a)\tend 0$ as $n\tend \infty$, if $n_1$ is large enough, we have
%% $\log d_{n_1}(\a)\geq \log r_{n_1}(\a)-\e/6$. Finally, by theorem
%% \ref{limitcremer}, if $n_1$ is large enough
%% $$\Phi_{n_1}(\a)+\log r_{n_1}(\a)\geq \Upsilon(\a)-\e/6.$$
%% Passing to the limit in the right side of equation (\ref{limitcremer1}) as
%% $\a'\in {\cal B}$  tends to $\a$, we obtain
%% $$(\forall \e>0)\quad
%% \liminf_{\a'\to \a,~\a'\in {\cal B}} \Phi(\a')+\log r(\a')
%% \geq \Upsilon(\a)-\e.$$

\section{Proof of inequality (\ref{liminf}) when
the Pérez-Marco sum diverges.\label{cremerparabo}}

In this section, we assume that $\a$ is a \Cremer number such that
$$\sup_n \beta_{n-1}\log \frac{1}{\a_n}=\infty.$$
To deal with this case, we will have to combine techniques of
parabolic explosion and techniques of renormalization.

Note that if $\beta_{n-1}\log 1/\a_n \leq C<\infty$ for all $n\geq
0$, then $\beta_{n-1}\log\log (e/\a_n)\leq
\beta_{n-1}\log(1+C/\beta_{n-1})$ decreases exponentially fast,
and $\a$ belongs to the set of \Cremer numbers studied in
section~\ref{loglogsec}.

\subsection{Parabolic explosion.}

The techniques of parabolic explosion are used to have a precise
control on the position of some periodic points of $P_{\a'}$ for
$\a'$ close to $\a$. The maps $P_{\a'}$, for $\a'$ real, are injective
on $B(0,1/2)$. We let $F_{\a'} \in S(\a')$ be
the lift of $P_{\a'}$ via $Z\mapsto \frac{1}{2}e^{2i\pi Z}$.
Let us recall that we called a periodic point of a map $F$ that commutes
with $T$, a point $Z$ such that $F^q(Z) = p$ for integers $q\in\N^*$
and $p\in \Z$ ($p$ and $q$ need not be coprime). Then $q$ is called
the period, and $p/q$ the rotation number.

\begin{lemma}\label{posper}
There exists a constant $B_\a>0$ such that for all Brjuno number $\a'$
sufficiently close to $\a$ and all integer $n\geq 2$,
\begin{itemize}
\item[a)] if $\frac{1}{2\pi}\Phi_n(\a')-B_\a>0$, then
$F_{\a'}$ has a periodic point with
period $\leq q_n$ and imaginary part $h'_0\geq
\frac{1}{2\pi}\Phi_n(\a')-B_\a$
\item[b)] in the upper half-plane
$$\left\{Z\in \C~|~\im(Z)\geq \frac{1}{2\pi}
\Phi_{n-1}(\a')+B_\a\right\}$$ 
the periodic points $Z$ of $F_{\a'}$ of period less
than or equal to $q_n$ come from ${\cal
C}_{p_n/q_n}(\a')$ (in the sense that $\frac{1}{2}e^{2i\pi Z} \in
{\cal C}_{p_n/q_n}(\a')$). 
\end{itemize}
\end{lemma}

\proof
For $n\geq 2$ and for $\a'\in \R\setminus \Q$, let us define
$X_n^*(\a')=X_n(\a')\setminus {\cal
  C}_{p_n/q_n}(\a')$  $r_n^*(\a')=\rad(X_n^*(\a'))$,
$d_n(\a')=d(0,X_n(\a'))$ and $d_n^*(\a')=d(0,X_n^*(\a'))$.

By proposition \ref{decprop} (since $q_2\geq 2$), we have for
$\a'$ close enough to $\a$,
$$\Phi_n(\a')+\log r_n(\a') \leq \Phi_2(\a')+\log r_2(\a') +C\sum_{k=3}^n
\frac{\log q_k}{q_k}.$$
As $\a'\tend \a$, the right hand term is bounded independently of $n$.
So, there exists a constant $C_\a$ such that for all $n\geq 2$ and all
$\a'\in {\cal B}$ sufficiently close to $\a$,
$$\Phi_n(\a')+\log d_n(\a')\leq \Phi_n(\a')+\log r_n(\a')\leq 2\pi C_\a.$$
Thus, if $\a'$ is sufficiently close to $\a$, $F_{\a'}$ has a
periodic point with imaginary part $h'_0\geq \ds \frac{1}{2\pi}
\Phi_n(\a')-C_\a-\frac{\log 2}{2\pi}$ when the right hand is positive.
This proves part a).

By lemma \ref{multiple}, in $B(\a',1/2q_n^3)$, the only cycle of
period less than or equal to $q_n$ that does not move
holomorphically is the cycle ${\cal C}_{p_n/q_n}(\a')$.
So, as in lemma \ref{equarat}, for all $n\geq 2$, we have
$$\Phi(\a')+\log r(\a') \leq \Phi_{n-1}(\a') + \log r_n^*(\a') +
(C-1)\sum_{k\geq n} \frac{\log q'_k}{q'_k},$$ where $C$ is the
constant provided by lemma \ref{lemmatec}. By
inequality (\ref{yoccoz}), $\Phi(\a')+\log r(\a')$ is universally
bounded from below. So, there exists a constant $C'$ such that for
all $n\geq 2$ and all $\a'$ sufficiently close to $\a$,
$$\Phi_{n-1}(\a') + \log r_n^*(\a') \geq -C'.$$

Finally, we claim that there exists a constant $C'_\a$ such that for
all $n\geq 2$ and all $\a'$
sufficiently close to $\a$, we have
$$\log d_n^*(\a')\geq \log r_n^*(\a') -C'_\a.$$
Part b) follows easily.
%% Indeed, let $z\in X_n^*(\a')$ be a point  that realizes the distance
%% $d(0,X_n^*(\a'))$ and set $w=P_{\a'}(z)$.
%% Then,
%% $$r_n^*(\a') \leq \rad(\C\setminus \{z,w\}) = d_n^*(\a') \cdot
%% \rad\left(\C\setminus \left\{1,\frac{w}{z}\right\}\right).$$ As
%% $\a'$ tends to $\a$, $w/z$ remains in a compact subset of
%% $\C\setminus \{0,1\}$ and so, $\rad(\C\setminus \{1,w/z\})$ is
%% bounded.
To prove the claim, let $\rho'=e^{2i\pi\a'}$ and $\rho=e^{2i\pi\a}$.
Let $n_0$ be such that $d_{n_0}^*(\a)<|\rho-1|/4$ (this is possible
since $\a$ is a \Cremer number). For $\a'$ close enough to $\a$,
$d_{n_0}^*(\a')<|\rho'-1|/2$.
For each fixed value of $n<n_0$, $\log d_n^*(\a') -\log
r_n^*(\a') \tend \log d_n^*(\a) -\log r_n^*(\a)$ when $\a' \tend
\a$. For $n\geq n_0$, let $z\in X_n^*(\a')$ be a point  that realizes
the distance $d_n^*(\a')$ and set $w=P_{\a'}(z) =\rho' z +
z^2$. Then, $|z|=d_n^*(\a')\leq d_{n_0}^*(\a')<|\rho'-1|/2$ and
$$r_n^*(\a') \leq \rad(\C\setminus \{z,w\}) = d_n^*(\a') \cdot
\rad\left(\C\setminus \left\{1,w/z\right\}\right).$$ As
$\a'$ tends to $\a$, $w/z = \rho'+z$ remains in a compact subset of
$\C\setminus \{1\}$ and so, $\rad(\C\setminus \{1,w/z\})$ is
bounded.
\qedprop

\subsection{Renormalization.}

Let us now fix $\d\in ]0,1/10[$. For $n\geq 0$, we will define a
sequence of heights $t'_n$ and a sequence of maps $F_{\a',n}\in
S(\a'_n)$ as in section \ref{secyoc}.

Let us set $$C'=2\pi(B_\a+4C_\d+32\d),$$ where $B_\a$ is the
constant in lemma \ref{posper}.

Now, let us choose $n_0$ so that $\beta_{n_0-1}\log
1/\a_{n_0}>4C'$ (this is possible
because $\sup \beta_{n-1}\log 1/\a_n=\infty$).
If $\a'$ is sufficiently close to $\a$, we
have
$$\beta'_{n_0-1}\log 1/\a'_{n_0}>4C'.$$

By proposition \ref{yocprop}, we can take
$$t'_0=\frac{1}{2\pi}\log \frac{1}{\a'_0}+C_\d\quad\ldots\quad
t'_{n_0-1}=\frac{1}{2\pi}\log \frac{1}{\a'_{n_0-1}}+C_\d.$$
By lemma \ref{posper} part~a), $F_{\a'}$ has a periodic point $Z'_0$ with
period $\leq q_{n_0}$ satisfying
$\im(Z'_0)=h'_0\geq \ds \frac{1}{2\pi}\Phi_{n_0}(\a')-B_\a$.
Note that $$\frac{1}{2\pi}\Phi_{n_0}(\a')-B_\a \geq
\frac{1}{2\pi}\Phi_{n_0-1}(\a')+\frac{4C'}{2\pi}-B_\a \geq
\frac{1}{2\pi}\Phi_{n_0-1}(\a')+B_\a.$$ 
By lemma \ref{posper} part~b), this periodic point comes from ${\cal
C}_{p_{n_0}/q_{n_0}} (\a')$, and thus has rotation number $p_{n_0}/q_{n_0}$.
By proposition~\ref{percor},
$F_{\a',n_0}$
has a fixed point $Z'_{n_0}$ with $\im(Z'_{n_0})=h'_{n_0}$ satisfying
$$\beta'_{n_0-1}h'_{n_0} + \sum_{n=0}^{n_0-1} \beta'_{n-1} t'_n +32\d
> h'_0$$ (see inequality (\ref{heightdn})). So,
$$h'_{n_0}>\frac{1}{2\pi}\log
\frac{1}{\a'_{n_0}}-\frac{B_\a+4C_\d+32\d}{\beta'_{n_0-1}}>\frac{3}{4}
\cdot\frac{1}{2\pi}\log \frac{1}{\a'_{n_0}}.$$ If $Z\neq Z'_{n_0}$
is another fixed point of $F_{\a',n_0}$, then proposition~\ref{percor} and
lemma \ref{posper} imply that
$$\beta'_{n_0-1}\im(Z) + \sum_{n=0}^{n_0-1} \beta'_{n-1}t'_n <
\frac{1}{2\pi}\sum_{n=0}^{n_0-1}
\beta'_{n-1} \log \frac{1}{\a'_n} +B_\a.$$
Thus,
$$\im(Z)<\frac{B_\a}{\beta'_{n_0-1}}<\frac{1}{4}\cdot\frac{1}{2\pi}
\log \frac{1}{\a'_{n_0}}.$$
So, there is a gap of height greater than $\ds \frac{1}{2}\cdot \frac{1}{2\pi}
  \log \frac{1}{a'_{n_0}}$ that separates the fixed point $Z'_{n_0}$ of
  $F_{\a',n_0}$ from the other fixed points of $F_{\a',n_0}$.
According to the second remark after proposition \ref{pmprop}, we
can therefore take
$$t'_{n_0} = h'_{n_0} + 1+ C_\d.$$
Finally, for $n\geq n_0+1$, we can take
$$t'_n=\frac{1}{2\pi}\log \frac{1}{\a'_n}+C_\d.$$

As in the previous section, proposition~\ref{size} we have
\begin{eqnarray*}
\log \frac{1}{2r(\a')}  & \leq &
2\pi\left(\sum_{n=0}^\infty \beta'_{n-1}t'_n + 32\d\right)\\
& \leq &  2\pi\left(\sum_{n=0}^{n_0-1} \beta'_{n-1}t'_n +
\beta'_{n_0-1}h'_{n_0}\right)+
\sum_{n=n_0+1}^{\infty} \beta'_{n-1}\log \frac{1}{\a'_n}\\
&& +2\pi(\beta'_{n_0-1}(4C_\d+1)+32\d)
\\
& \leq & 2\pi
h'_0+\Phi(\a')-\Phi_{n_0}(\a')+2\pi(\beta'_{n_0-1}(4C_\d+1)+32\d).
\end{eqnarray*}
Note that $2\pi h'_0 \leq - \log (2d_{n_0}(\a'))$ where
$d_{n_0}(\a')=d(0,X_{n_0}(\a'))$. So, reordering the terms and
simplifying by $\log 2$, we get
$$\Phi(\a')+\log r(\a')\geq \Phi_{n_0}(\a')+\log
d_{n_0}(\a')-2\pi(\beta'_{n_0-1}(4C_\d+1)+32\d).$$ We can now
conclude as in section \ref{loglogsec}.

\appendix

\section{Extracts from \cite{bc2}}\label{app_A}

The following proposition is proposition 10 from \cite{bc2}.

\begin{proposition}\label{prop_BC2_10}
Assume $U,V\subset \C$ are two hyperbolic domains containing $0$
and $\chi:U\to V$ is a holomorphic map fixing $0$. Let $S$ be a
finite subset of $U$ avoiding $0$, such that $\chi(S)$ avoids $0$.
Then,
$$\frac{\rad(V\setminus \chi(S))}{\rad(V)}\leq \frac{\rad(U\setminus
  S)}{\rad(U)}.$$
\end{proposition}

Given an integer $q\geq 1$, set
$$\U_q = \left\{e^{2i\pi k/q}~\big|~k=0,\ldots,q-1\right\}.$$
The following proposition is proposition 12 from \cite{bc2}.

\begin{proposition}\label{prop_BC2_12}
There exists a constant $C>0$ such that for $q\geq 2$ and $r<1$, we have
$$\log\rad(\D\setminus r\U_q) \leq \log r + \frac{C}{q}.$$
one can take $C=\log 4 + 2\log(1+\sqrt{2})$.
\end{proposition}

Let $V_\lambda$ be hyperbolic subdomains of $\C$ which contain $0$
and move holomorphically with respect to $\lambda\in \D$. The
following proposition is proposition 13 from \cite{bc2}.

\begin{proposition}\label{prop_BC2_13}
There exists a family of simply connected open sets $\wt{V}_\lambda$ and of
universal coverings $\pi_\lambda : \wt{V}_\lambda \to V_\lambda$ such
that $\wt{V_0} = \D$, the set
$$\wt{\cal V} = \big\{ (\lambda,z) \in \D \times\C\, \big| z \in
\wt{V}_\lambda \big\}$$ is open,
and $\Pi : (\lambda,z) \in \wt{\cal V} \mapsto \pi_\lambda(z)$ is
analytic.\\
For all $\lambda \in \D$,
$$\wt{V}_\lambda \subset B(0,\rho) \text{ with }
\log \rho = \frac{2\, \log 4}{\ds 1+ |\lambda|^{-1}}.$$
\end{proposition}

\section*{Aknowledgments.}

We would like to thank J.C. Yoccoz for several fruitful discussions
and suggestions.

\newcounter{nom}{\setcounter{nom}{1}}

\end{document}